\documentclass[10pt]{article}
\usepackage{makeidx}
\usepackage[utf8]{inputenc}
\usepackage[T1]{fontenc}
\usepackage{todonotes}
\usepackage{stmaryrd}
\usepackage{latexsym,amsfonts,amsthm,amsmath, amscd,amssymb,color}
\usepackage{float}
\usepackage[colorlinks=true,citecolor=blue,linkcolor=black,]{hyperref}
\usepackage[all]{xy}

\usepackage{geometry}
\usepackage{csquotes}
\MakeOuterQuote{"}

\theoremstyle{plain}
\newtheorem{lemma}{Lemma}[section]
\newtheorem{theorem}[lemma]{Theorem}
\newtheorem{proposition}[lemma]{Proposition}
\newtheorem{corollary}[lemma]{Corollary}
\theoremstyle{definition}
\newtheorem{definition}[lemma]{Definition}
\newtheorem{example}[lemma]{Example}
\newtheorem{remark}[lemma]{Remark}
\newtheorem*{definition*}{Definition}

\theoremstyle{remark}

\newtheorem{convention}[lemma]{Convention}
\newtheorem{notation}[lemma]{Notation}
\newtheorem{question}[lemma]{Question}

\usepackage{cite}

\usepackage{titling}

\title{On Nash resolution of (singular) Lie algebroids}

\author{Ruben Louis%\\ email \href{mailto:louisruben96@yahoo.com}{louisruben96@yahoo.com}
\thanks{Department of Mathematics, Jilin University, Changchun 130012, Jilin, China}~\thanks{Institut f\"ur Mathematik, Georg-August-Universit\"at G\"ottingen, G\"ottingen, Germany.}}

\begin{document}
\maketitle

\vspace{1em}

\begin{abstract}
Any Lie algebroid $A$ admits a Nash-type blow-up $\mathrm{Nash}(A)$ that sits in a nice short exact sequence of Lie algebroids $0\rightarrow K\rightarrow \mathrm{Nash}(A)\rightarrow \mathcal{D}\rightarrow 0$ with $K$ a Lie algebra bundle and $\mathcal{D}$ a Lie algebroid whose anchor map is injective on an open dense subset. The base variety is a blow-up determined by the singular foliation of $A$ considered recently in Omar Mohsen's \cite{OmarMohsen}. We provide concrete examples. Moreover, we extend the construction to singular subalgebroids in the sense of Androulidakis-Zambon \cite{MZ}.\end{abstract}

\scalebox{0.8}{\textbf{Keywords}: Lie algebroids, singular foliations, singularities.}\\
\scalebox{0.8}{\textbf{MSC classes}: 53-XX, 53D17, 32S45.}\\\scalebox{0.8}{\textbf{ORCID}: 0000-0002-3582-7994}
\tableofcontents
%\newpage 
%\begin{center}\textbf{Statements and Declarations}\end{center}I declare that this paper has been composed solely by myself and that it has not been submitted, in whole or in part, in any previous journal. This paper is part of a post-doc fellowship at Jilin University and Göttingen University. I acknowledge important discussions with Jun Jiang and Camille Laurent-Gengoux at earlystages of the project. Finally, I would like to thank Yunhe Sheng for his technical assistance at the Mathematics Laboratory of Jilin University in the realization of this paper\vspace*{4em}\noindent\hfill%\begin{center}\end{center}\newpage
\section*{Introduction}
We address the following question,
\begin{question}
    Is the slogan “A Lie algebroid splits in a purely manifold part and a purely Lie algebra part” valid?
\end{question}
There is of course a naive answer: given a Lie algebroid $(A, [\cdot\,,\cdot]_A, \rho)$ over $M$ with anchor $\rho\colon \Gamma(A)\rightarrow \mathfrak X(M)$, there is an exact sequence
 $$\xymatrix{ \ker(\rho) \ar@{^(->}[r]& \Gamma(A)\ar@{->>}[r]&\mathrm{Im}(\rho).}
$$
Above, $\ker(\rho)\subseteq \Gamma (A)$ is a totally intransitive object, i.e., a purely Lie algebra object. Let us be more precise. Let $\rho_m\colon A_m\to T_mM$ be the anchor map of $A$ at a point $m\in M$. Consider the collection of Lie algebras $\{\mathrm{Sker}\rho_m\}_{m\in M}$ where $u\in \mathrm{Sker}\rho_m\subseteq \ker \rho_m$ if and only if there exists a local section $\Tilde u$ of $A$ through $u$ such that $\rho (\Tilde{u})=0$. The vector space $\mathrm{Sker}\rho_m$ is a Lie ideal of the Lie algebra $\ker \rho_m$. Hence, it is a Lie algebra and its Lie bracket $\{\cdot\,,\cdot\}_{\mathrm{Sker}\rho_m}$ is given by
 $$\{a(m), b(m)\}_{\mathrm{Sker}\rho_m}:=[a,b]_A(m),\qquad a,b\in \ker\rho$$ hence depends only on the point $m$. On the base manifold $M$, the image $\mathrm{Im}(\rho)$ is a singular foliation in the sense of \cite{AS,Cerveau, Debord,LLS}, so it is entirely “within” the manifold $M$. Hence, it is a purely manifold object. But this answer \emph{is not} satisfactory, since neither $\ker(\rho)$ nor $\mathrm{Im}(\rho) $ are sections of a vector bundle in general. However, we will argue that the slogan above is valid, provided that we first apply some canonical blow-up procedure invented by Omar Moshen \cite{OmarMohsen}, and in alignment with the conceptual framework of the Nash blow-up of an affine variety\cite{Nobile}.

Let us explain the results of this paper. An important class of Lie algebroids consists of those whose anchor maps are injective on an open dense subset. Those Lie algebroids have been extensively studied by Claire Debord \cite{debord2000local, Debord}. In this paper, we shall call them \emph{Debord Lie algebroids}. {Let us briefly discuss some key features of Debord Lie algebroids. By definition, the anchor map of a Debord Lie algebroid is injective at the level of sections, i.e., $\Gamma(A) \stackrel{\rho}{\hookrightarrow} \mathfrak{X}(M)$. As a result, the induced singular foliation $\mathcal{F} = \rho(\Gamma(A)) \subset \mathfrak{X}(M)$ is a finitely generated projective $C^\infty(M)$-module, since $\Gamma(A) \simeq \mathcal{F}$. Conversely, given any finitely generated singular foliation $\mathcal{F} \subseteq \mathfrak{X}(M)$ that is projective as a $C^\infty(M)$-module, there exists a vector bundle $A$ such that $\Gamma(A) \simeq \mathcal{F}$, by the Serre–Swan Theorem \cite{SwanRichardG, MoryeArchanaS}. This vector bundle carries a natural Lie algebroid structure whose anchor map $A \to TM$ is injective on an open dense subset. Therefore, the singular foliation associated with a Debord Lie algebroid is completely determined by the image of its anchor map. In this sense, it is a purely “manifold object”. One of their key properties is that they are always integrable to smooth Lie groupoids, in contrast to arbitrary Lie algebroids, which may not be integrable \cite{Debord, Crainic-Fernandes}. Debord Lie algebroids are naturally the class of Lie algebroids closest to regular foliations, that is, to vector subbundles $F\subseteq TM$ such that $\left [\Gamma(F), \Gamma(F)\right]\subset \Gamma(F)$.}

We present a method for desingularizing a Lie algebroid $A$ by finding a pair $(N,p)$ consists of a space $N$ and a map $p\colon N\to M$ such that\begin{enumerate}
    \item $p\colon N\to M$ is onto and proper;
    \item  the restriction $p|_{M_\mathrm{reg}}\colon p^{-1}(M_{\mathrm{reg}})\to M_{\mathrm{reg}}$ to regular points $M_\mathrm{reg}$ of $A$  is one-to-one
    \end{enumerate} on which the pullback $p^!A$ of $A$  admits a Lie algebroid structure $\left(\mathrm{Nash}(A),[ \cdot \,, \cdot ]_{\mathrm{Nash}(A)}, \Hat{\rho}\right)$ whose quotient is a Debord Lie algebroid. Also, $\mathrm{Nash}(A)\to A$ is a Lie algebroid morphism. In fact, by construction, $N$ comes equipped with an action of $A$ along $p\colon N\rightarrow M$ in the sense of \cite{T.Mokri,kosmann2002differential,Mackenzie}. The fundamental assertion presented in Proposition \ref{prop:main} of this paper is that we can choose $p\colon N\rightarrow M$ such that there exists a proper Lie algebroid ideal $K\rightarrow N$ of $p^!A$ (i.e., with $\Gamma(K)\subseteq\ker (\Hat{\rho})$) whose fibers coincide with those of $ \ker(\Hat{\rho})$ on the open dense subset of regular points. As a consequence, the quotient $\mathrm{Nash}(A)/K$ is a Lie algebroid with injective anchor map on an open dense subset. One of the most important invariants of this construction is the basic singular foliation of $A$ that is, the singular foliation given by the subsheaf of Lie-Rinehart algebras $\mathcal{F}:=\rho_A(\Gamma(A))\subseteq \mathfrak X(M)$. In other words, $N$ does not depend on the Lie algebroid $A$ but only on its basic singular foliation. Importantly, this result remains valid even when the underlying structure $(A\rightarrow  M , [\cdot\,, \cdot ]_A, \rho_A)$ is merely an almost Lie algebroid, i.e., without the requirement that $[ \cdot \,, \cdot ]_A$ adheres to the Jacobi identity: the anchor map still has to be a morphism of brackets.

%Our construction is inspired by the work of O. Mohsen \cite{OmarMohsen} on singular foliations also by other related works on desingularisation of Lie groupoids, see e.g., \cite{DebordClaire2017BcfL, NistorVictor2019DoLg}.
In some cases, our blow-up space $N$ is precisely the usual blow-up along the singular part of the singular foliation $\mathcal F = \rho_A(\Gamma(A))$ induced by the Lie algebroid. We must admit that our $N$ is often not a smooth manifold. However, our Lie algebroids structures on $N$ admit leaves that are smooth manifolds. But if the closed subset of singular points $M_{\mathrm{sing}}$ is a smooth manifold and $\ker\rho$ is finitely generated, then $N$ is often a smooth manifold, see \cite{Ali-Sinan} and \cite[\S 2.2.1]{Ruben3} for a more precise result.

The construction of $N$ relies on a very classical method which is due to the mathematician J. Nash \cite{Nobile}, used mainly in algebraic geometry for desingularisation. The idea of applying this method to singular foliations is not completely new. As far as I know,  A. Sert\"oz Sinan was the first who adapted it on foliations in his Ph.D. \cite{Ali-Sinan} to compute the Baum-Bott residues of singular holomorphic foliations, see also \cite{JPTS}. In fact,  A. Sert\"oz Sinan did more, he applied Nash in the realm of coherent subsheaves of locally free sheaves and generalized the works of Nobile \cite{Nobile}. The idea of Nash blow-up $\mathrm{Nbl}_\mathcal F(M)$ of a singular foliated manifold $(M, \mathcal{F})$ in the sense of \cite{AS} was used by O. Mohsen \cite{OmarMohsen, mohsen2022tangent} and consists in replacing every singular point by the limiting positions of the tangent spaces of the leaves in neighborhoods of regular points. In this paper, we do not apply the Nash method directly on the tangent space of our singular foliation $\mathcal{F}\subseteq \mathfrak X(M)$ but on the kernel of the anchor map $ \rho\colon A\rightarrow TM$. The choice of considering the kernel allows us to generalize the Nash construction to Lie $\infty$-algebroid that are universal in the sense of \cite{LLS,CLRL}. This is discussed by the author in \cite{louis2023universalhigherliealgebras,Ruben3}. The latter provides a class of applications of this method on solvable singular foliations, the author described a sequence of Nash blow-ups associated with a geometric resolution of a singular foliation in the sense of \cite{LLS}. %The latter is not necessarily locally finitely generated, in contrast with the case of coherent sheaf assumption of \cite{Ali-Sinan}. 

These constructions enable us, in passing, to efficiently compute several instances of blow-up of a singular foliation.

Also, we show in Proposition \ref{prop:Poisson} that the space of all limits of the tangent spaces of the regular symplectic leaves of a Poisson manifold is isomorphic to the space of all the limits of the kernel of the associated contangent algebroid on regular points. This justifies the name “Nash blow-up”. In the Poisson case, we give an example where the Poisson bivector field degenerates under the blow-up procedure, see Example \ref{ex:spheres}. However, all its Hamiltonian vector fields remain well-defined after this procedure.

Moreover, we claim that our method can be applied to many "Lie algebroids" type structures,  such as "almost Lie algebroids" \cite{MarcelaPaul},  "left-symmetric algebroids" in the sense of \cite{Jiefeng} or "Hom-Lie algebroids" \cite{LAURENTGENGOUX201369,ZHANG}, and so on. And one could easily replace, with no additional assumptions, the (Lie) bracket of our main theorem with the corresponding operation.

The structure of the paper is the following. In Section \ref{sec:1}, we recall some basics on singular subalgebroids and singular foliations. In Section \ref{sec:2}, we introduce the Nash construction and state the main results of the paper with their proofs. In Section \ref{sec:3}, we end the paper with some examples from singular foliations and Poisson geometry.\\

\noindent
\textbf{Acknowledgements}. This work was carried out as part of a joint postdoctoral fellowship between Jilin University and University of Göttingen. I sincerely thank the National Natural Science Foundation of China for awarding me the Research Fund for International Young Scientists. I acknowledge insightful discussions with Jun Jiang and Camille Laurent-Gengoux during the early stages of the project. I would like to thank Yunhe Sheng and Chenchang Zhu for their valuable technical assistance at the Departments of Mathematics at Jilin University and the University of Göttingen, which contributed significantly to the completion of this work. Finally, I thank the anonymous referee for his/her careful reading of the manuscript and for the constructive comments that helped improve the presentation.

\section{Definitions}\label{sec:1}
Throughout the paper, $M$ stands for a smooth/analytic/complex manifold or an algebraic variety over $\mathbb K\in \{\mathbb R, \mathbb C\}$. We denote by $\mathcal{O}_M$ the algebra of smooth/holomorphic or polynomial functions on $M$ depending on the context. Also, $\mathbb K$ stands for $\mathbb R$ or $\mathbb C$.\\

Let $(E, \left[\cdot\,,\cdot\right]_E, \rho_E)$ be a Lie algebroid on  a manifold $M$.
\subsection{Singular subalgebroids}
The following concept was first introduced in \cite{MZ}.
\begin{definition}\label{def:Sing-subal}
    A \emph{singular subalgebroid} of $E$ is an $\mathcal{O}_M$-submodule $\mathcal{F}\subseteq \Gamma (E)$ which is
    \begin{enumerate}
        \item \emph{involutive}, i.e., $[\mathcal{F}, \mathcal{F}]_E\subseteq \mathcal{F}$;
        \item \emph{locally finitely generated}\footnote{In the complex or real analytic case this is automatically satisfied \cite{Jean-Claude}.}, i.e.,  every point $m\in M$ admits an open neighborhood $\mathcal{U}$  on which there exists $\xi_1,\dots, \xi_k \in \mathcal{F}|_\mathcal{U}\subseteq \Gamma(E)|_\mathcal{U}$ such that for every $\mathcal V \subset \mathcal U$, the restrictions of $\xi_1, \dots, \xi_k$ to $\mathcal V$ generate $\mathcal{F}|_\mathcal V$ as a $\mathcal O_\mathcal V$-module.
    \end{enumerate}
\end{definition}

\begin{remark}
    Notice that Definition \ref{def:Sing-subal} means that a singular subalgebroid $\mathcal{F}\subseteq \Gamma(E)$ is a locally finitely generated Lie-Rinehart algebra \cite{HuebschmannJohannes} whose Lie bracket and anchor map are the restriction of the Lie bracket $[\cdot\,,\cdot]_E$ and the anchor map $\rho_E$ to $\mathcal{F}$. 
\end{remark}
\begin{example}[Singular foliations] One important class of singular subalgebroids for this paper is the class of singular foliations in the sense of \cite{Hermann,Debord,AS, LLL1}. A \emph{singular foliation} $\mathcal{F}$ on $M$ is a singular subalgebroid of $E=TM$. Here are some important features of singular foliations in the smooth/complex cases.
\begin{itemize}
     \item[-]Singular foliation admits leaves \cite{Hermann}: there exists a partition of $M$ into immersed submanifolds (possibly of non-constant dimension) called \emph{leaves} such that for all $m\in M$, the image of the evaluation map $\mathcal F \to T_m M  $  is the tangent space of the leaf through $m$.
     \item[-] \emph{Singular foliations are self-preserving}: the flow $\phi^X_t$ of vector fields $X\in\mathcal F $, whenever defined, preserves $\mathcal F$ \cite{Hermann,AS,GarmendiaAlfonso}, i.e., 
         $\forall\, m\in M, \exists\, \varepsilon >0$  and  a neighborhood $\mathcal{U}$ of $m$ such that $\forall t\in ]-\varepsilon,\varepsilon[$\begin{equation*} {(\phi^X_t|_\mathcal{U})}_*(\mathcal F|_\mathcal{U})=\mathcal F|_{\phi_t^X(\mathcal{U})}.
     \end{equation*}
\end{itemize}
\end{example}

\begin{example}[Debord singular subalgebroids] Here is a special class of singular subalgebroids. A singular subalgebroid $\mathcal{F}$ of $E$ is said to be \emph{Debord} or \emph{projective} if there is a vector bundle morphism  $A\stackrel{\rho}{\rightarrow} E$ covering the identity of $M$ such that $\rho(\Gamma(A))=\mathcal{F}$ and  which is injective on an open dense subset. {In other words, $\mathcal{F}\simeq \Gamma(A)$ is projective as an $\mathcal{O}_M$-submodule of $\Gamma(E)$.} In that case, there is a Lie bracket on $A$ such that $\rho$ is a Lie algebroid morphism. Therefore,  $\mathcal{F}$ arises as the image of a Lie algebroid whose anchor map is injective on an open dense subset.
\end{example}

For more  examples, see \cite{MZ}.

\subsubsection{$E$-valued almost Lie algebroids and singular subalgebroids}
Most of our method will be valid for almost Lie algebroids. Let us recall some definitions.

\begin{definition}[\cite{MarcelaPaul,LLS}]\label{def:LA}Let $\mathcal{F} \subseteq \Gamma(E)$ be a singular subalgebroid. An \emph{$E$-valued almost Lie algebroid  $(A, \left[\cdot\,,\cdot\right]_A,\rho_A)$  over $\mathcal{F}$} is a vector bundle morphism  $A\stackrel{\rho_A}{\rightarrow} E$ (referred as the \emph{anchor map}) such that $\rho_A(\Gamma(A))=\mathcal{F }$; together with a skew-symmetric bilinear map  $\left[\cdot\,,\cdot\right]_A\colon \Gamma(A)\wedge\Gamma(A)\rightarrow \Gamma(A)$ (referred to as \emph{the bracket}) satisfying for all $a,b\in \Gamma(A)$ and $f\in \mathcal{O}_M$
\begin{enumerate}

        \item $\rho_A(\left[a,b\right]_A)=[\rho_A(a),\rho_A(b)]_E$;
        \item $\left[a,fb\right]_A=\rho_E\circ\rho_A(a)[f]b+ f \left[a,b\right]_A$.
    \end{enumerate}
$\mathcal{F}$ shall be called the \emph{basic} singular subalgebroid of $(A, \left[\cdot\,,\cdot\right]_A,\rho_A)$. If $E=TM$ we simply say "almost Lie algebroid".
\end{definition}
If $\left[\cdot\,,\cdot\right]_A$ satisfies the Jacobi identity, i.e., for  $a,b, c\in \Gamma(A)$\begin{equation}
    [[a,b]_A,c]_A+ [[c,a]_A,b]_A+ [[b,c]_A,a]_A=0,
\end{equation}then  $A$ is a {Lie algebroid} over $M$ with anchor map $\rho_E\circ \rho_A$.
\begin{remark}Every singular subalgebroid is the image of an $E$-valued almost Lie algebroid on any open subset where it is finitely generated. For instance,  a singular subalgebroid is finitely generated if and only is it is the image of some $E$-valued almost Lie algebroid, see \cite[Proposition 3.8]{LLS}
\end{remark}
\begin{proposition}\label{prop:equivalence}
For any two  almost algebroids $(A_1, \left[\cdot\,,\cdot\right]_{A_1},\rho_{A_1})$ and $(A_2, \left[\cdot\,,\cdot\right]_{A_2},\rho_{A_2})$ over a singular subalgebroid $\mathcal{F}\subseteq \Gamma(E)$,  there exist morphisms of vector bundles \begin{equation}
    \xymatrix{A_1\ar@<2pt>[r]^\Phi&A_2\ar@<2pt>[l]^\Psi}
\end{equation}
such that
\begin{enumerate}
    \item $\rho_{A_2}\circ \Phi=\rho_{A_1}$ and $\rho_{A_1}\circ \Psi=\rho_{A_2}$.
      \item $\rho_{A_2}(\Phi\circ\Psi -\mathrm{id})=0$ and $\rho_{A_1}(\Psi\circ\Phi -\mathrm{id})=0$.
\end{enumerate}
In particular, $[ \Phi(a), \Phi(b)]_{A_2}- \Phi([a,b]_{A_1} ) \in \ker \rho_{A_2}$ for all $a,b \in \Gamma(A_1)$. Also, we have $[ \Psi(a), \Psi(b)]_{A_1}- \Psi([a,b]_{A_2} ) \in \ker \rho_{A_1}$ for all $a,b \in \Gamma(A_2)$.
\end{proposition}
\begin{proof}
For any pair of $E$-valued almost Lie algebroids $ (A_1, [\cdot\,,\cdot]_{A_1}, \rho_{A_1})$ and $({A_2},[\cdot\,,\cdot]_{A_2},\rho_{A_2})$ over $\mathcal{F}$, there exists (by projectivity of the module of sections of vector bundles) vector bundle morphisms \begin{equation*}
    \xymatrix{(A_1, [\cdot\,,\cdot]_{A_1}, \rho_{A_1})\ar@<2pt>[r]^\Phi&(A_2, [\cdot\,,\cdot]_{A_2}, \rho_{A_2})\ar@<2pt>[l]^\Psi}
\end{equation*}
that satisfies item 1. The other properties are automatically verified.
\end{proof}
Let us give a crucial result on linear vector fields on $E$-valued almost Lie algebroids. These techniques are used in \cite{zbMATH07105909}. Recall that for $q\colon A\rightarrow M$ a vector bundle over $M$, a \emph{linear vector field} on $A$ is a couple $(\xi, X)\in \mathfrak X(A)\times \mathfrak X(M)$ such that $$\xymatrix{A\ar[r]^\xi \ar[d]_q&TA\ar[d]^{Tq}\\M\ar[r]^X&TM}$$ is a  morphism of vector bundles. Equivalently,  \begin{equation}
     \xi[\Gamma(A^*)]\subseteq \Gamma(A^*)\;\;\text{and}\; \;\xi[q^*\mathcal O_M]\subseteq q^*\mathcal O_M.
 \end{equation}

 For more details and properties on linear vector fields on vector bundles, see \cite{Mackenzie}. 

\begin{proposition}\label{prop:lift:1}
    Let $(A, \left[\cdot\,,\cdot\right]_A,\rho_A)$ be an $E$-valued almost Lie algebroid  over a singular subalgebroid  $\mathcal{F}\subseteq \Gamma(E)$. Let  $e\in \mathcal{F}=\rho(\Gamma(A))$ and $X=\rho_E(e)$,  there exist a linear vector field $\xi_e$ on $q\colon A\rightarrow M$ and a linear vector field $\nu_X$ on $q'\colon E\to M$ such that
 \begin{enumerate}
\item $\xi_e$ is $q$-related to $X$ and $\nu_X$ is $q'$-related to $X$,
\item their flows at time $t$ are defined if and only if the flow of $X$ at time $t$ is defined,
 \item the flows $\phi^{\xi_e}_t,\;\phi^{\nu_X}_t$ whenever defined are isomorphism of vector bundles that make the following diagram commutes, 
\begin{equation}\label{eq:diag1}
\xymatrix{A\ar[r]^{\phi_t^{\xi_e}}\ar[d]_{\rho_A}&A\ar[d]^{\rho_A}\\ E\ar[d]_{\rho_E}\ar[r]^{\phi_t^{\nu_X}}&E\ar[d]^{\rho_E}\\TM\ar[r]^{T\phi_t^X} &TM.}
\end{equation}
\end{enumerate}
\end{proposition}
\begin{proof}
Let $e\in\mathcal F$ and $X=\rho_E(e)$. Choose $a\in \Gamma(A)$ a section of $q\colon A\to M$ such that $\rho_A(a)=e$. The Lie derivatives $\mathcal{L}_a,\,\mathcal{L}_e$ defined  by 
\begin{enumerate}
    \item \begin{align*}
    \mathcal{L}_a[q^*f]:&=q^*(X[f]),\;\forall\; f\in \mathcal O_M,\\\langle\mathcal{L}_a[\alpha], b\rangle&= X[\langle \alpha, b\rangle]-\langle \alpha, [a,b]_{A}\rangle, \quad \forall \alpha\in \Gamma(A^*),\; b\in \Gamma(A);
 \end{align*}
 \item \begin{align*}
    \mathcal{L}_e[q'^*f]:&=q'^*(X[f]),\;\forall\; f\in \mathcal O_M,\\\langle\mathcal{L}_e[\beta], e'\rangle&= X[\langle \beta, e'\rangle]-\langle \beta, [e,e']_{E}\rangle, \quad \forall \beta\in \Gamma(E^*),\; e'\in \Gamma(E)
 \end{align*}
\end{enumerate}
define linear vector fields $\xi_e\in \mathfrak X(A),\; \nu_X\in\mathfrak X(E)$ which are respectively $q$-related to $X$  and $q'$-related to $X$. Notice that these linear vector fields depend, respectively, on the  brackets $[\cdot\,,\cdot]_A,\; [\cdot\,,\cdot]_E$  and the vector field $X$. This proves item 1. %An argument to prove item 2 and 3 is the following: the vector field $(\xi_e+\nu_X, X)$ is a linear vector field on $A\oplus E$ which is tangent to the graph of $\rho_A\colon A\to E$$$ \left\{ (u, \rho_A(u) ) | u \in A \right\} \subset A \oplus E. $$  Hence, the diagram \eqref{eq:diag1} commutes. For a complete proof of item 2 and 3 see  \cite{LLL}, Proposition 2.2.11. p. 108.
{An argument to prove item 2 and 3 is to show that $(\xi_e+\nu_X, X)$ is a linear vector field on $A\oplus E$ that is tangent to the graph of $\rho_A\colon A\to E$
 $$ \left\{ (u, \rho_A(u) ) | u \in A \right\} \subset A \oplus E. $$  Hence, the diagram \eqref{eq:diag1} commutes. Indeed, for $\beta\in \Gamma(E^*)$ and $f\in \mathcal{O}_M$ and $b\in \Gamma(A)$, we have \begin{align*}
     \langle \mathcal{L}_a[\beta\circ \rho_A],b\rangle&=X\langle\beta\circ \rho_A, b\rangle-\langle \beta\circ \rho_A,[a,b]_A\rangle\\&=X\langle\beta, \rho_A( b)\rangle-\langle \beta, \rho_A([a,b]_A)\rangle\\&=X\langle\beta, \rho_A( b)\rangle-\langle \beta, [\rho_A(a),\rho_A(b)]_E\rangle\\&=\langle\mathcal{L}_e[\beta],\rho_A(b)\rangle= \langle\mathcal{L}_e[\beta]\circ\rho_A, b\rangle.
 \end{align*}
 This means that $\mathcal{L}_a[\beta\circ \rho_A]=\mathcal{L}_e[\beta]\circ\rho_A$ for all $\beta\in \Gamma(E^*)$. Likewise, \[\mathcal{L}_a[{q'}^*f\circ \rho_A]=\mathcal{L}_e[{q'}^*f]\circ\rho_A\] for all $f\in \mathcal{O}_M$. This implies that the vector fields $\xi_a\in \mathfrak X(A)$ and $\nu_X\in \mathfrak X(E)$ are $\rho_A$-related. Therefore, the vector field $\xi_a+ \nu_X$ is tangent to the graph of $\rho_A\colon A\to E$, since the vector fields $\xi_a$ and $\xi_a+\nu_X$ are related via the embedding\begin{align*}
     &A\longrightarrow A\oplus E\\&u\,\longmapsto (u, \rho_A(u)).
 \end{align*}
 This completes the proof.}
\end{proof}
\begin{remark}
    When $E=TM$, we have $\phi_t^{\nu_X}=T\phi_t^X$ in Diagram \eqref{eq:diag1}.
\end{remark}

\subsubsection{Isotropy Lie algebras of a singular subalgebroid}
We start with the following definition

\begin{definition}[Strong-kernel]
    Let $\eta\colon A\rightarrow E$ be a vector bundle morphism covering the identity map. The \emph{strong-kernel} $\mathrm{Sker}\,\eta_m$ of $\eta$ at $m\in M$ is the subvector space of $\ker\eta_m$ defined as $u\in A_m$ such that there is an open neighborhood $\mathcal{U}$ of $m$ and a local section $\Tilde{u}\in\Gamma(A)|_{\mathcal{U}}$ through $u$ so that $\eta(\Tilde{u})=0$.
\end{definition}
Here are some easy properties of the strong-kernel
\begin{proposition}\cite{LLL2,AS}
 Let $\eta\colon A\rightarrow E$ be a vector bundle morphism covering the identity map.

 \begin{enumerate}
     \item The map $m\mapsto \dim (\mathrm{Sker}\,\eta_m)$ is lower semi-continuous and $m\mapsto \dim (\ker\eta_m)$ is upper semi-continuous.

     \item  $\mathrm{Sker}\,\eta_m=\ker\eta_m$ if and only if there is an open neighborhood $\mathcal{U}$ of $m$ such that $\mathrm{Sker}\,\eta_{m'}=\ker\eta_{m'}$ for all $m'\in \mathcal{U}$.
     \item The open dense subset where $m\mapsto \dim (\ker\eta_m)$ and $m\mapsto \dim(\mathrm{Im}\,\eta_m)$ are
locally constant is the set of points $m\in M$ such that $\mathrm{Sker}\,\eta_m=\ker\eta_m$.
 \end{enumerate}
\end{proposition}

The notion that we now introduce generalizes the original notion of isotropy Lie algebras for singular foliations in \cite{AS}.

Let $\mathcal
F\subseteq \Gamma(E)$ be a singular subalgebroid. We may assume that $\mathcal{F}$ is finitely generated, so that it is the image of an $E$-valued almost Lie algebroid $(A, \left[\cdot\,,\cdot\right]_A, \rho)$. For $m\in M$, the skew-symmetric bracket $\left[\cdot\,,\cdot\right]_A$ induces a bilinear map $$\left[\cdot\,,\cdot\right]_{A,m}\colon \wedge^2 \ker\rho_m\rightarrow \ker\rho_m$$ defined for all $u, v\in \ker \rho_m$ by $\left[u\,,v\right]_{A,m}:= \left[\widetilde u\,,\widetilde v\right]_A$ where $\widetilde u, \widetilde v$ are respectively local sections through $u$ and $v$ respectively in an open neighborhood of $m$. It is easily checked that this pointwise bracket is well-defined, that is $\left[\cdot\,,\cdot\right]_A$ depends on the $1$-jet at $m$ of the sections of $A$. Furthermore, the strong kernel $\mathrm{Sker}\,\rho_m$ at $m$ satisfies
 $ [{\mathrm{Sker}}\,\rho_m,{\ker}\rho_m]_{A,m} \subset {\mathrm{Sker}}\,\rho_m$. This implies that the skew-symmetric bilinear map $ [ \cdot, \cdot  ]_{A, m} $ goes to the quotient to a bilinear map,
  \begin{equation}\label{eq:isLie} \{\cdot, \cdot\}_m \colon  
  \wedge^2 \frac{\ker\rho_{m}}{\mathrm{Sker}\,\rho_m} \longrightarrow
   \frac{\ker\rho_{m}}{\mathrm{Sker}\,\rho_m}
  \end{equation}
which turns out to be a Lie bracket \cite{LLL2}.

\begin{definition}
 The \emph{isotropy Lie algebra of $\mathcal{F}$ at $m$} is the quotient
 $$ \mathfrak g_m (\mathcal F) =
 \frac{\ker\rho_{m}}{\mathrm{Sker}\,\rho_m}$$
 together with the Lie bracket $\{\cdot\,,\cdot\}_m$.
\end{definition}
The following proposition is a direct consequence of Proposition \ref{prop:equivalence}
\begin{proposition}
    The isotropy Lie algebra $\mathfrak g_m (\mathcal F)$ of a singular subalgebroid $\mathcal{F}$ at $m$ is independent of the choice of an $E$-valued almost Lie algebroid over $\mathcal F$.
\end{proposition}

\begin{proof}
This is a direct consequence of Proposition \ref{prop:equivalence}.
\end{proof}
\begin{remark}[Androulidakis-Skandalis isotropy Lie algebras] Let  $\mathcal{\mathcal{F}}\subseteq \Gamma(E)$ be a singular subalgebroid and  $m\in M$. Denote by $\mathcal{I}_m$ be the ideal of vanishing functions at $m$ and $\mathcal{F}(m)$ the Lie algebra of sections of $\mathcal{F}$ vanishing at $m$ which contains the ideal $\mathcal{I}_m\mathcal{F}$. Assume that $\mathcal{F}$ is the basic singular subalgebroid of some $E$-valued almost Lie algebroid $(A, [\cdot\;,\cdot]_A, \rho)$. There is a natural Lie algebra isomorphism $\mathfrak g_m(\mathcal{F})\simeq \frac{\mathcal{F}(m)}{\mathcal{I}_m\mathcal{F}}$,  see \cite{LLL2}. When $\mathcal{F}\subseteq \mathfrak X(M)$ is a singular foliation, $\mathfrak g_m(\mathcal{F})$ is the isotropy Lie algebra of $\mathcal{F}$ at $m$ in the sense of \cite{AS}. 
\end{remark}
\section{Nash blow-up of an $E$-valued almost Lie algebroid}\label{sec:2}
\subsection{Grassmann bundle}

%Let $d, k \in  \mathbb N$ be such that $0\leq k\leq d$ and $\mathbb K= \mathbb R, \mathbb C$. The set  $\mathrm{Grass}_{k}(\mathbb K^d)$  of all $k$-dimensional vector subspaces of $\mathbb K^d$ is a complex compact manifold of dimension $k(d-k)$, called  the \emph{Grassmanian of $k$-planes in  $\mathbb K^d$} see e.g., \cite{Bredon}. The case $k=1$ corresponds to the projective space $\mathbb P^{d-1}(\mathbb K)$.Let us spell out the manifold structure of $\mathrm{Grass}_{k}(\mathbb K^N)$ for $\mathbb K=\mathbb C$. The groups ${U}(d)$ of invertible $\mathbb C$-linear maps, and of unitary linear transformations, acts transitively on $\mathrm{Grass}_{k}(\mathbb C^d)$ by $ g \cdot V = g(V) $for all $V\in\mathrm{Grass}_{k}(\mathbb C^d)$ and $g \in {U}(d)$. This action is well-defined since $g$ is invertible. Therefore,  \begin{eqnarray*}\mathrm{Grass}_{k}(\mathbb C^d) & \stackrel{set}{\simeq}& \frac{{ U}(d)}{ U(k)\times  U(d-k)} \end{eqnarray*} Above the subgroups $U(k)\times  U(d-k)$   is the  stabilizer of $\mathbb C^k \times \{0\}^{d-k} \subset \mathbb C^d$ in $U(d)$.  Since both groups in the first line above are compact Lie groups, the first description equips $\mathrm{Grass}_{k}(\mathbb C^N)$ with a structure of compact manifold. If one replaces $\mathbb C$ by $\mathbb R $ and $U(d)$ by $O(d)$,  then $\mathrm{Grass}_{k}(\mathbb R^d)$ is a compact smooth manifold.

\begin{convention}
For our current purpose, we consider the Grassmannian of all sub-spaces of \underline{co}-dimension $k\in\mathbb N$ in $\mathbb C^N $ rather than the set of subspaces of a given dimension. It is convenient to denote this manifold by $\mathrm{Grass}_{-k}(\mathbb C^N)$ (notice the use of a minus sign). In other words, we set $\mathrm{Grass}_{-k}(\mathbb C^N) := \mathrm{Grass}_{N-k}(\mathbb C^N)$.
\end{convention}

Now, we consider this construction pointwise on a given vector bundle over a manifold. The next definition makes sense in the smooth or complex cases without adaptation.

\begin{definition}
Let $A\rightarrow M$ be a vector bundle of rank $d$ over a manifold $M$. Let $k\leq \mathrm{rk}(A)=d$. The disjoint union:
  $$\mathrm{Grass}_{-k}({A}):= \coprod_{x\in M} \mathrm{Grass}_{-k}(A_{x})  $$
  comes equipped with a natural manifold structure. Also,
 \begin{equation}\label{eq:Pi}
     P\colon \mathrm{Grass}_{-k}(A)\longrightarrow M
 \end{equation}
 is a fibration with typical fiber $\mathrm{Grass}_{-k}(\mathbb K^d) $. It is called \emph{$-k$-th Grassmann bundle of~$A$}.
\end{definition}
 
 %\begin{remark}For every open subset $\mathcal U\subset M$ on which $A$ is trivial, $P^{-1}(\mathcal U)\simeq \mathcal{U}\times\mathrm{Grass}_{-k}(\mathbb{K}^d).$  An \emph{adapted chart} for $\mathrm{Grass}_{-k}({A})\longrightarrow M$ around a point $x\in M$ is a then defined as a set of local coordinates of the form $(P^*x_1,\ldots P^*x_n)$ and $(z_{i,j})_{i \in [1:k] , j \in [1: n-k]}$ where $(x_1,\ldots,x_n)$ are local coordinates on $M$ and $(z_1,\ldots, z_{r(d-k)})$ are functions which are standard affine coordinates on an open subset of each fiber of $P$. \end{remark}
 
 The following elementary proposition will be used in the subsequent sections. Again, it is valid in both smooth and complex cases. 
\begin{proposition}
    Let $A\longrightarrow M$ be a vector bundle over $M$. A linear vector field $(\xi, X)$ on $A$ with base vector field $X\in \mathfrak X(M)$ induces a vector field on the Grassmann bundle $P\colon \mathrm{Grass}_{-k}(A)\longrightarrow M$ which is $P$-related to $X$.
\end{proposition}

\noindent
\textbf{The tautological bundle}.
There is a natural and important exact sequence of vector bundles over $\mathrm{Grass}_{-k}(A)$, for $A \to M$ a vector bundle and $k$ an integer $ \leq {\mathrm{rk}}(A)$.

\begin{enumerate}
\item Let $ P^!A$ be the pullback  of the vector bundle $A \to M$ to $ \mathrm{Grass}_{-k}(A)$.
\item The \emph{tautological subbundle} $\tau_A^{-k}\longrightarrow\mathrm{Grass}_{-k}(A)$ is the vector subbundle of $P^! A $ whose fiber at a point  $V_x\in\mathrm{Grass}_{-k}(A)$ is precisely the vector space $V_x$, seen as a subspace of $\left(P^! A\right)_{V_x}=A_x$  (the $ -k$ being a reminder that it is of corank $k$ in $P^! A$).
\item  We denote by $Q_A^k$ and call \emph{tautological quotient bundle} the bundle over $\mathrm{Grass}_{-k}(A)$ obtained by taking the quotient of the first bundle by the second one, that is, $Q_A^k:= P^! A / \tau_A^{-k}$ (the "$k$" being a reminder that this quotient has rank $k$).
\end{enumerate}
These three vector bundles fit into the short exact sequence
\begin{equation}\label{eq:shortsequence-vb1}
    \xymatrix{
      \tau_A^{-k}\ar@{^(->}[r] \ar[d]&\ar[d] P^! A \ar@{->>}[r] &  Q^k_A \ar[d] \\\mathrm{Grass}_{-k}(A)\ar@{=}[r]& \mathrm{Grass}_{-k}(A)\ar@{=}[r] & \mathrm{Grass}_{-k}(A)
 }
\end{equation} 
\subsection{The blow-up space $\mathrm{Nbl}_\mathcal{F}(M)$: Nash construction} \label{sec:Nbl-construction}
Omar Mohsen \cite{OmarMohsen} developed a blow-up method for an action Lie algebroid associated with a Lie algebra $\mathfrak{g}$ acting on a manifold $M$. Additionally, he employed a comparable approach to define the blow-up space of any arbitrary singular foliation, without using additional underlying structures. In this section, we exploit his concept to build a blow-up space using anchored bundles. This facilitates the application of this methodology to any object resembling a Lie algebroid. As said in the Introduction, the construction is valid for arbitrary vector bundle morphisms, which is a particular case of A. Sert\"oz Sinan \cite{Ali-Sinan} applied in the realm
of coherent subsheaves of a locally free sheaf.  This point of view allows us to apply our result to a larger class than the one of singular foliations, namely, singular subalgebroids in the sense of Androulidakis-Zambon \cite{MZ}.

Let $A, E$ be vector bundles over $M$. Consider a morphism of vector bundles $$\xymatrix{A\ar[r]^{\rho}\ar[d] &E\ar[d]\\ M\ar@{=}[r]&M}$$
Denote by $\mathcal{F}:=\rho(\Gamma(A))\subseteq \Gamma(E)$ the subsheaf of $\mathcal{O}_M$-submodules (after sheafification).
\begin{convention}
    We denote by $M_{\mathrm{reg}} \subseteq  M$ the open dense subset where $x\mapsto \dim (\ker\rho_x)$ and $x\mapsto \dim(\mathrm{Im}\,\rho_x)$ are
locally constant. The points of $M_{\mathrm{reg}}$ are called \emph{regular} and the points in $M\setminus M_{\mathrm{reg}}$ are called \emph{singular}. In the complex case, the dimension of $\mathrm{Im}\,\rho_x\subseteq E_x$ is the same for all points $x\in M_{\mathrm{reg}}$. In the smooth case, we shall now assume that $\rho$ is of constant rank on $M_{\mathrm{reg}}$. This happens in particular when $\mathcal{F}$ admits a geometric resolution of finite length \cite{LLS}. Let $r$ be the common dimension of $\mathrm{Im}\,\rho_x\subseteq E_x$ on $ M_{\mathrm{reg}}$. The restriction $\mathcal{F}|_{M_\mathrm{reg}}\subseteq E|_{M_\mathrm{reg}}$ is a vector subbundle of rank $r$, called the \emph{regular part} of $\mathcal{F}$.
\end{convention}
 %we assume that the regular leaves of $\mathcal{F}$ have the same dimension $r$. 
\begin{remark}
    When $E=TM$ and $\mathcal{F}\subseteq \mathfrak X(M)$  is a singular foliation, $M_{\mathrm{reg}}$ is the union of the regular leaves of $\mathcal{F}$. More generally, when $E$ is a Lie algebroid and $\mathcal{F}\subseteq \Gamma(E)$ is a singular subalgebroid, $M_{\mathrm{reg}}$ is the set of points $m\in M$ whose isotropy Lie algebra vanishes, i.e., such that $\mathfrak g_m(\mathcal{F})=0$. Moreover, the restriction $\mathcal{F}|_{M_\mathrm{reg}}\subseteq E|_{M_\mathrm{reg}}$ is a Lie subalgebroid of $E|_{M_\mathrm{reg}}$.
\end{remark}

\noindent
\textbf{Construction of the blow-up space}. Notice that for every $x\in M_{\mathrm{reg}}$
\begin{enumerate}
\item $\mathrm{Im}\,\rho_x$ is a point of the Grassmannian  $\mathrm{Grass}_{r}(E_x)$, and
\item $\ker \rho_x$ is a point of $\mathrm{Grass}_{-r}(A_x)$.
\end{enumerate}
Consider now the Grassmann bundle  $\mathrm{Grass}_{r}(E)$, and  the Grassmann bundle $\mathrm{Grass}_{-r}(A)$. We denote by $P$, in both cases, their natural projections onto $M$. Consider the two natural sections of $P$ defined on $M_{\mathrm{reg}}$ by: 
\begin{enumerate}
\item  $\sigma_{\mathrm{Im}}\colon M_{\mathrm{reg}}\longrightarrow \mathrm{Grass}_{r}(E),\, x\longmapsto \mathrm{Im}\,\rho_x$, and
\item  $\sigma_{\ker}\colon M_{\mathrm{reg}}\longrightarrow \mathrm{Grass}_{-r}(A),\, x\longmapsto \ker\rho_x$.
\end{enumerate}
Then we define 
\begin{enumerate}
\item  $\mathrm{Nbl}^{\mathrm{Im}}_\mathcal{F}(M):= \overline{\sigma_{\mathrm{Im}}(M_{\mathrm{reg}})}$ to be the closure of the image of the section $\sigma_{\mathrm{Im}}$ in $\mathrm{Grass}_{r}(E)$, and
\item $\mathrm{Nbl}^{\ker}_\mathcal{F}(M):= \overline{\sigma_{\ker}(M_{\mathrm{reg}})}$ to be the closure of the image of the section $\sigma_{\ker}$ in $\mathrm{Grass}_{-r}(A)$.
\end{enumerate}

We denote by $p$ the restriction of $P$ to both $\mathrm{Nbl}^{\mathrm{Im}}_\mathcal{F}(M)$ and $\mathrm{Nbl}^{\ker}_\mathcal{F}(M)$. 

\begin{remark}
Notice that in \cite{Ruben3}, a sequence, depending on $ n \in \mathbb N$, of "blow-ups" is constructed, which for $n=0$ gives back $\mathrm{Nbl}^{\mathrm{Im}}_\mathcal{F}(M)$, and for $n=1$ gives $\mathrm{Nbl}^{\ker}_\mathcal{F}(M)$. 
\end{remark}

In Section \ref{application:poisson}, we have a special class of examples where $\mathrm{Nbl}^{\ker}_\mathcal{F}(M)\simeq \mathrm{Nbl}^{\mathrm{Im}}_\mathcal{F}(M)$, namely Poisson structures.
\begin{convention}
    From now on, we will only consider $\mathrm{Nbl}^{\ker}_\mathcal{F}(M)$, which we simply denote by $\mathrm{Nbl}_\mathcal{F}(M)$ and shall be called \emph{Nash blow-up} of $M$ along $\mathcal{F}$. This notation, which only mentions $(M, \mathcal{F})$, is justified in Proposition \ref{prop:well-defined} below.
\end{convention}

\begin{remark}
\label{rmk:limitesGrassmann}
Notice that\begin{enumerate}
    \item  Intuitively, for $x\in M$, $p^{-1}(x)= \mathrm{Nbl}_\mathcal{F}(M)\cap P^{-1}(x)$ is the set of all possible limits of the subspaces $\ker \rho_{y}$ when $y\in M_{\mathrm{reg}}$ converges to $x$.
More precisely, for any $x\in M$, there is an open neighborhood $\mathcal U\subset M$ of $x$ such that $\mathrm{Grass}_{-r}(A)\simeq \mathcal{U}\times \mathrm{Grass}_{-r}(\mathbb{K}^{\mathrm{rk}(A)})$. By construction,   \begin{equation*}
        p^{-1}(x)=\left\{V\subset A_{_x}\;\middle|\; \exists\, (x_n)\in M_{\mathrm{reg}}^\mathbb N,\, \text{such that},\,\;\ker \rho_{x_n} \underset{n \to +\infty}{\longrightarrow} V \text{ as }  x_n\underset{n \to +\infty}{\longrightarrow}x\right\}.
        \end{equation*}
        \item For any open subset $\mathcal{U}$ of $M$, $\sigma_{\ker}(\mathcal{U}\cap M_{\mathrm{reg}})$ is dense in $p^{-1}(\mathcal{U})$.
\end{enumerate}
\end{remark}

In the following proposition, we gather some easy but important features of that construction. 
\begin{proposition}\label{prop:properties}
Let $\rho\colon A\to E$ be a vector bundle morphism. In the smooth case, we assume that all regular leaves have the same dimension $r$. The map $p\colon \mathrm{Nbl}_{\mathcal {F}}(M)\rightarrow M$ (referred as the \emph{blow-up map}) satisfies the following properties:
 \begin{enumerate}
\item 
$p$ is proper and surjective. In particular, for each point $x\in M$, the fiber $p^{-1}(x)$ is a non-empty compact set.
\item {For every $x\in M$ and $V\in p^{-1}(x)$,  $\mathrm{Sker}\,\rho_x\subseteq V\subseteq \ker\rho_x$.} 
\item For every $x\in M$ and $V\in p^{-1}(x)$, the image $ [V]$ of $V$ in the quotient vector space %isotropy Lie algebra at $x$, i.e., 
$$ \mathfrak g_x (\mathcal F) =
 \frac{\ker\rho_{x}}{\mathrm{Sker}\,\rho_x}$$ is a vector space %Lie subalgebra 
 of codimension $r-\dim (\mathrm{Im}\,\rho_x)$. %$r-\dim(L_x)$, where $\dim (L_x)$ is the dimension of the leaf through $x$.
 Also, $V \to [V ] $ is an injective map from $p^{-1}(x) $ to ${\mathrm{Grass}}_{-(r - \dim (\mathrm{Im}\,\rho_x))} (\mathfrak g_x (\mathcal F)  )  $.
\item  For every $x\in M_{\mathrm{reg}}$,  $p^{-1}(x)=\{\ker\rho_x\}$ is reduced to a single point in $\mathrm{Grass}_{-r}(A)$. Also, $p^{-1}(M_{\mathrm{reg}})$
is a smooth manifold and the restriction $p\colon p^{-1}(M_{\mathrm{reg}})\longrightarrow M_{\mathrm{reg}}$ is invertible\footnote{Invertible here means: diffeomorphism, in the smooth case, bi-holomorphism, in the complex case.}.
\end{enumerate}

\end{proposition}
\begin{proof}
    The proof follows the same lines as Proposition 2.4 in \cite{Ruben3} for $n=1$. %But \cite{Ruben3} assumed existence of a geometric resolution for $\mathcal{F}$ in the sense of \cite{LLS}. However, for $n=1$ in \cite{Ruben3} the latter assumption is not needed, in other words, $\ker(\rho)$ does not need to be finitely generated. Indeed, if $A\rightarrow TM$ is the beginning of a geometric resolution of $\mathcal{F}$ that is, if there exists vector bundle morphism $$E\stackrel{\mathrm{d}}{\rightarrow} A\stackrel{\rho}{\rightarrow} TM$$ such that $\rho\circ \mathrm{d}=0$ and $\mathrm{Im}(\mathrm{d})=\ker\rho$ at the level of the sheaf of sections, then we have $\mathrm{Sker}(\rho_x)=\mathrm{Im}(\mathrm{d}_x)$ for all $x\in M$, in particular $\mathfrak g_x(\mathcal{F})=\frac{\ker\rho_x}{\mathrm{Im}(\mathrm{d}_x)}$.
\end{proof}

\begin{remark}
    Despite the extremely pleasant properties, the blow-up space $\mathrm{Nbl}_\mathcal{F}(M)$ has a major problem: it is not a submanifold in general. However, since $\mathrm{Nbl}_\mathcal{F}(M)\subset \mathrm{Grass}_{-r}(A)$  is defined as a closed subset, it inherits a “good” algebra of smooth functions and can be seen as a nice differential space, see Section \ref{sec:differential-space}. In addition, we show that it admits in the case we are interested in a partition into smooth manifolds which form  in fact the leaves of the pullback  foliation\footnote{We make sense of singular foliation on $\mathrm{Nbl}_\mathcal{F}(M)$ in Section \ref{sec:differential-space}.}of $\mathcal{F}$ see Section \ref{sec:lifting}.
\end{remark}

{We have the following smoothness result, which is a direct consequence of \cite[Corollary 1.2]{Ali-Sinan} and has been further refined in \cite[\S 2.1.2]{Ruben3}.
\begin{proposition}\label{prop:smoothness} Let $k\in \mathbb N$ be the dimension of $\{\ker\rho_m\}_{m\in M}$ on $M_{\mathrm{reg}}$.
    If \begin{enumerate}
        \item $\ker\rho\subset \Gamma(A)$ is a finitely generated $\mathcal{O}_M$-submodule and the singular locus $\Sigma=M\setminus M_{\mathrm{reg}}$ is a smooth submanifold of $M$, \item and the ideal generated by the  $k\times k$ minors of sections of $A$ that generate $\ker \rho$ on an open subset is $\mathcal{I}_\Sigma^n$ for some $n\in \mathbb N$ (here $\mathcal{I}_\Sigma$ is the ideal of vanishing functions on $\Sigma$), 
    \end{enumerate}then $\mathrm{Nbl}_\mathcal{F}(M)$ is the usual blow-up of $M$ along $\Sigma=M\setminus M_{\mathrm{reg}}$. Hence, $\mathrm{Nbl}_\mathcal{F}(M)$ is a smooth manifold.
\end{proposition}}
%\begin{remark}The first assumption of Proposition \ref{prop:smoothness} is fulfilled in particular if $\mathcal{F}$ admits a geometric resolution in the sense of \cite{LLS}.  For a more general statement than Proposition \ref{prop:smoothness} see \cite{Ruben3}. % If $\mathcal F$ admits a geometric resolution of lenght $n+1$ and $Z=M\setminus M_{\mathrm{reg}}$ is smooth, then $\Tidle M_n$ is smooth.\end{remark}

Until now, we had defined the Nash blow-up for any vector bundle morphism $\rho\colon A\to E$. Now, let us add an $E$-valued almost Lie algebroid bracket on $A$ so that the image $\mathcal{F}=\rho(\Gamma(A))$ is a singular subalgebroid of $E$.

%\begin{convention}In the sequel, we take $E=TM$ and $\mathcal{F}$ a singular foliation on $M$, but the constructions remain the same without any adaptation in the singular subalgebroid case.  \end{convention}
The proof of the following proposition follows the same lines as in Proposition 2.4 in \cite{Ruben3}.
\begin{proposition}\label{prop:almost}
    Let $\rho\colon A\to E$ be a vector bundle morphism. 
We assume that $(A\rightarrow M, \rho)$ is equipped with a skew-symmetric bracket that makes it an $E$-valued almost Lie algebroid with anchor map $\rho$. In the smooth case, we assume that all regular leaves have the same dimension $r$. 
 \begin{enumerate}

\item For every $x\in M$ and $V\in p^{-1}(x)$, $V$ is stable under the skew-symmetric bracket\footnote{If $A$ is  a Lie algebroid, then $\ker\rho_m$ is a Lie algebra and $V\subseteq \ker\rho_m$ is a Lie subalgebra of $\ker\rho_m$.} of $\ker\rho_m$.
\item For every $x\in M$ and $V\in p^{-1}(x)$, the image $ [V]$ of $V$ in the isotropy Lie algebra $ \mathfrak g_x (\mathcal F)$ at $x$ is a Lie subalgebra 
 of codimension $r - \dim (\mathrm{Im}\,\rho_x)$%, where $\dim (L_x)$ is the dimension of the leaf through $x$
 .
 Also,  $V \to [V ] $ is an injective map from $p^{-1}(x) $ to ${\mathrm{Grass}}_{-(r - \dim (\mathrm{Im}\,\rho_x)) } (\mathfrak g_x (\mathcal F)  )  $.

\end{enumerate}
\end{proposition}

\begin{remark}
    Let $(A, [\cdot\,,\cdot]_A, \rho)$ be an $E$-valued almost Lie algebroid. The Jacobiator $J\colon \wedge^3A\rightarrow A$ (the failure of having Jacobi identity) fits into the non-exact complex \begin{equation}
        \wedge^3A\stackrel{J}{\rightarrow} A \stackrel{\rho}{\rightarrow} E
    \end{equation}
   Here is a natural question which can be the object of another study: what kind of structure do we obtain if we apply the Nash construction to the vector bundle morphism $J$?
\end{remark}

Now we justify the notation $\mathrm{Nbl}_\mathcal F(M)$: it does not mention the anchor
map $\rho\colon A \rightarrow E$, because, indeed, it does not depend on the $E$-valued almost Lie algebroid structure. More precisely, 
\begin{proposition}
\label{prop:well-defined}
$\mathrm{Nbl}_\mathcal{F}(M)$ does not depend on the choice of an $E$-valued almost Lie algebroid over $\mathcal{F}$, i.e., the sets associated with two different anchored bundles are canonically isomorphic.
\end{proposition}
\begin{proof}
Let $\mathcal{F}\subseteq \Gamma(E)$ a singular subalgebroid, and let $A, A'$ be two $E$-valued almost Lie algebroids over $\mathcal{F}$. Let  \textbf{ \begin{equation}
    \xymatrix{A\ar@<2pt>[r]^\Phi&A'\ar@<2pt>[l]^\Psi}
\end{equation}} as in Proposition \ref{prop:equivalence}. We denote temporarily by  $p_A\colon \mathrm{Nbl}_\mathcal{F}(M)\rightarrow M$ also by $p_{A'}\colon \mathrm{Nbl}_\mathcal{F}(M)\rightarrow M$) the induced projection map w.r.t to $A$ (resp. to $A'$). 
\begin{enumerate}
    \item Let $x\in M$. It is easily checked for all $V_x\in p_A^{-1}(x) $ and $V'_x\in p_{A'}^{-1}(x) $ that $\dim [{V}_x]=\dim [{V}_x']$, see \cite{Ruben3}. 
    \item Let $x\in M$ and $k=\dim (\mathrm{Sker}\rho_x)$. Let $(u_i)_{i=1, \ldots, k}$ be local sections of $A$ in an open neighborhood $\mathcal{U}_x$ of $x$ such that $\rho(u_i)=0$ and $\mathrm{{Span}}(u_i(x))=\mathrm{Sker}\rho_x$. We can shrink $\mathcal{U}_x$ so that the map $y\mapsto \mathrm{Span}(u_i(y))\subseteq \mathrm{Sker}\rho_y$  be of constant rank on $\mathcal{U}_x$. Therefore, $ K_A(x):=\sqcup_{y\in \mathcal{U}_x}\mathrm{Span}(u_i(y))\subset \ker\rho|_{\mathcal{U}_x}$ is a vector subbundle of $A$ of rank $k$ on $\mathcal{U}_x$ whose fiber over $x$ is $\mathrm{Sker}\rho_x$. Likewise, there exists a vector subbundle $K_{A'}(x)$ of $A'$ on $\mathcal{U}_x$  whose fiber over $x$ is $\mathrm{Sker}\rho'_x$. Without loss of generality, one can assume that $\Phi_{|_y}(K_{A}(x))\subseteq K_{A'}(x)$ for all $y\in \mathcal{U}_x$. %Since $\Phi(\mathrm{Sker}\rho)\subset \mathrm{Sker}\rho'$, 
Thus, the map $\Phi_{|_y}$ goes to quotient 
$$\phi_y\colon \frac{\ker\rho_y}{K_{A}(x)}\longrightarrow  \frac{\ker\rho'_y}{K_{A'}(x)}$$
for all $y\in \mathcal{U}_x$ so that $\phi_x=[\Phi_{|_x}]$. Let $(x_n)_{n\in \mathbb N}$ be a sequence of regular points converging to $x$ such that the sequence $\ker\rho_{x_n}$ converges to $V_x$ in $\mathrm{Grass}_{-r}(A)$. This implies that the sequence $\ker\rho_{x_n}/K_{A}(x)$ converges to $[V_x]=V_x/K_{A}(x)=V_x/\mathrm{Sker}\rho_x$  in $\mathrm{Grass}_{-r}\left(\frac{A}{K_{A}(x)}\right)$.  Since $\phi_{x_n}\left(\frac{\ker\rho_{x_n}}{K_{A}(x)}\right)\subseteq \frac{\ker\rho'_{x_n}}{K_{A'}(x)}$, it follows that $\phi_x([V_x])\subseteq [V'_x]$ where $V'_x$ is the limit of (a subsequence of) $\ker\rho'_{x_n}$.  By item 1, $[V_x]$ and $[V'_x]$ have the same dimension, therefore $\phi_x([V_x])= [V'_x]$.  This defines the required map and completes the proof. 
\end{enumerate}
    \end{proof}

\subsubsection{Non-smoothness of $\mathrm{Nbl}_\mathcal{F}(M)$ and Nagano-Sussman}\label{sec:differential-space}

Although the blow-up space $\mathrm{Nbl}_\mathcal{F}(M)\subseteq\mathrm{Grass}_{-r}(A)$ may fail to be a smooth manifold (see Example \ref{ex:non-smooth}), we explain how it can inherit a sort of “smooth structure”, i.e., it has a good enough framework to define the notion of singular foliations on $\mathrm{Nbl}_\mathcal{F}(M)$ that admit smooth leaves. 

%Although $\mathrm{Nbl}_\mathcal{F}(M)\subseteq\mathrm{Grass}_{-r}(A)$ may fail to be a  smooth manifold, %we show in Section \ref{sec:lifting} that it can be endowed with a notion of singular foliation.In this sectionwe explain how the blow-up space $\mathrm{Nbl}(M, \mathcal{F})$ can inherit a sort of “smooth structure” known as “differential space”\cite{Sniatycki} which is a notion weaker than smooth manifold. This smooth structure on $\mathrm{Nbl}_\mathcal{F}(M)$ allows having a good framework to define the  notion of  singular foliations  on $\mathrm{Nbl}_\mathcal{F}(M)$.

Since our blow-up space $\mathrm{Nbl}_{\mathcal {F}}(M)$ is defined in particular as a closed subset of $\mathrm{Grass}_{-r}(A)$, it can be equipped with a good notion of “smooth” functions \begin{equation}\label{eq:nash-functions}
    C^\infty\left(\mathrm{Nbl}_\mathcal F(M)\right):=\left\{g_{|\mathrm{Nbl}_\mathcal{F}(M)},\,\text{for all}\, \, g\in C^\infty\left(\mathrm{Grass}_{-r}{(A)}\right)\right\}
.\end{equation}If $\mathrm{Nbl}_\mathcal{F}(M)$ is a smooth submanifold, then it is an embedded submanifold, since it is a closed subset. Moreover, the algebra $C^\infty\left(\mathrm{Nbl}_\mathcal{F}(M)\right)$ defined in \eqref{eq:nash-functions}  coincides with the algebra of smooth functions on $\mathrm{Nbl}_\mathcal{F}(M)$.\\ %Unless stated otherwise, we will always consider $\mathrm{Nbl}_{\mathcal {F}}(M)$ with its subcartesian differential space structure. This is sufficient for us to state our results and examples.

We introduce the following definitions, which are particular cases of a more general notion known as subcartesian differential spaces, see \cite{Sniatycki}. We will work in the context of smooth differential geometry, but this can be adapted easily in the context of holomorphic/analytic cases.
\begin{enumerate}
%\item Let $S$ be a closed subset of a manifold $M$ and let $\mathcal{O}_S$ be the algebra made of the restrictions of smooth functions of $M$ to $S$. Let  $\xi\in \mathrm{Der}(\mathcal{O}_S)$ be a derivation of $\mathcal{O}_S$. An \emph{integral curve} of $\xi$  through $s\in S$ is a smooth curve $\gamma \colon I\subseteq \mathbb R \longrightarrow S\subseteq M$ such that for every $g\in \mathcal{O}_S$ $$\frac{d}{dt}g\circ \gamma(t)= \xi[g]\circ \gamma(t), \; \text{for all}\; t\in I.$$ One of the important features of this concept is that derivations of $\mathcalO_S$ admit maximal integral curves  through a given point. This allows to construct the flows. We shall denote by $\mathrm{exp}(t\xi)_s$ the integral curve at $t\in I\subseteq\mathbb R$ of $\xi$ through $s\in S$. 
        \item %The map $\mathrm{exp}(t\xi)\colon S\to S$ might fail to be a local diffeomorphism as shown in Example 1 of \cite{Sniatycki}. Hence, vector fields on $S\subseteq M$ can not be only derivations. 
            Let $S$ be a closed subset of a manifold $N$.  A \emph{vector field} on $S\subseteq N$  is the restriction to $S$ of a vector fields $Z\in \mathfrak X(N)$  whose flow preserves $S$, i.e., $\phi^{Z}_t(S)\subseteq S$ whenever it makes sense. In that case, we shall say that such a \emph{$Z$ is tangent to $S$}. The set of vector fields on $S$ form a Lie algebra which we denote as in the usual case by $\mathfrak X(S)$.
        \item The \emph{tangent space} $T_sS$  of $S$ at $s\in  S$ is the evaluation at $s$ of the vector fields on $S$.
\end{enumerate}
%\begin{example}Let $S\subseteq M$ a closed subset. A vector field $Z\in \mathfrak X(M)$ whose flows  $\phi_t^Z$  preserves $S$ when it makes sense, restricts to a vector fields on $S$.\end{example}

\begin{remark}
    Notice that when $S$ is a submanifold, this notion of vector fields on $S$ agrees to the usual case.
\end{remark}

Now,  we recall a crucial theorem that allows us to correctly define singular foliations on a closed subset $S$ of a manifold $N$.

\begin{definition}\cite{LLL1, LLS}\label{def:leaves}\label{def:SF2} Let $S$ be a closed subset of $N$.
\begin{enumerate}
    \item  A \emph{singular foliation on $S$} is an involutive\footnote{Notice that $\mathcal{F}$ is generated by the restrictions to $S$ of vector fields on $N$, they are  required to be involutive only after restrictions to $S$.} locally finitely generated $\mathcal{O}_S$-submodule $\mathcal{F}\subseteq\mathfrak X(S)$.
    \item For $s\in S$, the \emph{leaf} of a singular foliation $\mathcal{F}$ on $S$ through $s$ is the set 
\begin{equation}
        L_s:=\left\{\phi_{t_1}^{Z_1}\circ \phi_{t_2}^{Z_2}\circ \cdots \circ \phi_{t_k}^{Z_k}(s),\; t_1, \ldots, t_k\in \mathbb R\; \right\}
    \end{equation}
 Above, $Z_1,\ldots,Z_k$ are vector fields on $N$  whose restrictions are in $ \mathcal{F}$. We implicitly assume that the flows are defined.
\end{enumerate}
\end{definition}
\begin{remark}
    From Definition \ref{def:leaves}, it is easily checked that being in the same leaf is an equivalence relation on $S$, therefore the leaves induce a partition of $S$.
\end{remark}

The notion of leaves of a singular foliation on $S\subseteq M$ is justified by the following theorem, which generalizes the Stefan-Sussmann theorem \cite{Stepan1,Stepan2}, which says the leaves are smooth manifolds.

\begin{theorem}\label{thm:Sussmann_SP}Let $\mathcal{F}$ be a singular foliation on a closed subset  $S\subseteq M$. The leaves $\mathcal F$ form a partition of $S$ into connected manifolds,  immersed as submanifolds of $M$. %Each leaf $L\subseteq S$ admits a suitable topology that comes with a unique manifold structure such that the inclusion map $L\hookrightarrow S \subseteq M$ is smooth.
\end{theorem}
The explanation of this result is based on a very strong theorem known as the Nagano-Sussmann theorem \cite{Nagano}. Let us first recall this theorem.

The following theorem, widely used in control theory, provides a very strong result regarding the smoothness of the orbits of a finite number of vector fields on a manifold  without any assumptions. 
       \begin{theorem}[Nagano–Sussmann] Let $\mathcal{V}\subseteq \mathfrak X(N)$ be a locally finitely generated $\mathcal{O}_N$-submodule of vector fields on a manifold $N$. For every $\ell\in N$, the set
       $$\left\{\phi_{t_1}^{Z_1}\circ \phi_{t_2}^{Z_2}\circ \cdots \circ \phi_{t_n}^{Z_n}(\ell),\; t_1, \ldots, t_n\in \mathbb R,\; Z_1,\ldots, Z_n\in \mathcal{V},\;n\in \mathbb N\right\}$$

     is a connected immersed submanifold of $N$.
       \end{theorem}

\begin{proof}[Proof (of Theorem \ref{thm:Sussmann_SP})]
   %This theorem is much simpler in our context, because the vector fields on $S$ are the restrictions of vector fields on $M$ to $S$.
   We can assume that $\mathcal{F}\subseteq \mathfrak X(S)$ is globally finitely generated. Let $\xi_1,\ldots, \xi_k$ be generators for $\mathcal F$. By definition, the $\xi_i$'s are the restrictions to $S$ of vector fields $Z_i$'s on $N$ whose flows $\phi^{Z_i}_t$ preserves $S$, i.e.,   $\phi^{Z_i}_t(S)\subseteq S$ where the flows are defined. By Nagano-Sussmann theorem, the orbits generated by the vector fields  $Z_1,\ldots, Z_k\in \mathfrak X(N)$  are immersed submanifolds of $N$. By assumption, the orbit that contains a point of $S$ is included in $S$. This completes the proof.
\end{proof}

Likewise, using this framework, the notion of (almost) Lie algebroids can be naturally defined on a closed subset $S$ of a manifold $N$.
\begin{definition}\label{def:LA2}
    \emph{An almost Lie algebroid on  $S$} is the triplet $(A|_S, [\cdot\,,\cdot]_{A, S}, \rho_{A, S})$ where \begin{enumerate}
        \item $A|_S$ is the restriction of a vector bundle $A\to N$ to $S$;
\item $\rho_{A,S}$ is the restriction to $S$ of a vector bundle morphism $\rho_A\colon A\to TN$ such that $\rho(a)|_S$ defines a vector field on $S$, for all $a\in \Gamma(A)$;
        \item $[\cdot\,,\cdot]_{A, S}$ is skew-symmetric $\mathbb R$-bilinear map on $ \Gamma(A|_S)$ obtained as the restriction of a skew-symmetric  $\mathbb K$-bilinear map $[\cdot\,,\cdot]_A \colon \Gamma(A)\times \Gamma(A)\to \Gamma(A)$ satisfying for all $a,b\in \Gamma(A)$
        
        \begin{enumerate}
            \item the usual Leibniz identity w.r.t $\rho_A$ : $\left[a,fb\right]_A=\rho_A(a)[f]b+ f \left[a,b\right]_A$,
            \item $\rho_A([a,b]_A)-[\rho_A(a), \rho_A(b)]\in \mathcal{I}_S\mathfrak X(N)$,
        \end{enumerate}
        where $\Gamma(A|S)$ are the restrictions of sections of $A$ to $S$ and $\mathcal{I}_S\subset \mathcal O_N$ is the ideal of vanishing functions on $S$.
    \end{enumerate}
    If $[[a,b]_A,c]_A+[[c,a]_A,b]_A+ [[b,c]_A,a]_A\in\mathcal{I}_S\Gamma(A)$, for all $a,b,c\in \Gamma(A)$, we say  $(A|_S, [\cdot\,,\cdot]_{A, S}, \rho_{A, S})$ is a \emph{Lie algebroid on  $S$}.
\end{definition}
\begin{remark}
    For an almost Lie algebroid on a closed subset $S$ as in Definition \ref{def:LA2}, the image $\mathcal{F}=\rho_{A,S}(\Gamma(A|S))$ is a singular foliation on $S$ in the sense of Definition \ref{def:SF2}.
\end{remark}
\begin{example}
    In this context, a \emph{Debord Lie algebroid} on a closed subset $S\subseteq N$ is a Lie algebroid $(A|_S, [\cdot\,,\cdot]_{A, S}, \rho_{A, S})$ such that $\rho_{A,S}\colon A|_S\to TS $ is injective on an open dense subset of $S$. Also, a singular foliation $\mathcal{F}\subseteq \mathfrak X(S)$ is said to be a \emph{Debord foliation}  if it is the image of a Debord Lie algebroid on $S$. In that case, $\mathcal{F}=\rho_{A,S}(\Gamma(A|_S))$ is a projective $\mathcal{O}_S$-module.
\end{example}

\subsection{The (almost) Lie algebroid structure $\mathrm{Nash}(A)\to \mathrm{Nbl}_\mathcal{F}(M)$}\label{sec:lifting}
In this section, we state the main results of the paper and discuss their proofs.

Recall that in Section \ref{sec:Nbl-construction}, we constructed the blow-up space $\mathrm{Nbl}_\mathcal{F}(M)$ out of a finitely generated singular subalgebroid $\mathcal{F}\subseteq \Gamma(E)$, with the mapping $$p\colon \mathrm{Nbl}_\mathcal{F}(M)\longrightarrow M,$$ satisfying the conditions outlined in Proposition and \ref{prop:properties} and Proposition \ref{prop:almost} through an $E$-valued almost Lie algebroid $A\to M$, yet independent of it. In this section, we study the pullback of $A\to M$ to $\mathrm{Nbl}_\mathcal{F}(M)$. All constructions remain applicable for $E$-valued almost Lie algebroids. Subsequently, we narrow our focus to the case of Lie algebroids.\\

\begin{notation}In this section, \begin{enumerate}
\item $(E, \left[\cdot\,,\cdot\right]_E, \rho_E)$ is a Lie algebroid on  $M$,
    \item $(A, \left[\cdot\,,\cdot\right]_A, \rho)$ stands for an (almost) Lie algebroid valued in $E$, and its basic singular subalgebroid is denoted by $\mathcal{F}=\rho(\Gamma(A))\subseteq \Gamma(E)$.
\end{enumerate} Recall that the regular part of $\mathcal{F}$ has constant dimension $r$. Also, $p^!A$ stands for the restriction to $\mathrm{Nbl}_\mathcal{F}(M)$  of the vector bundle $P^!A\to \mathrm{Grass}_{-r}(A)$. The  pullback of a section  $a\in \Gamma(A)$ to $\Gamma(p^!A)$ shall be denoted by $p^!a$. \end{notation}

\subsubsection{Main statements}
The main result of this paper establishes that any singular subalgebroid $\mathcal{F}\subseteq \Gamma(E)$ arising from an $E$-valued almost Lie algebroid $A\to M$ is desingularized into a Debord Lie algebroid on  the blow-up space $p\colon \mathrm{Nbl}_{\mathcal{F}}(M)\to M$. This Lie algebroid embeds into an exact sequence, wherein the last two terms are almost Lie algebroids. The latter becomes a sequence of Lie algebroids if $A\to M$ is a Lie algebroid. %whose last term is Lie algebroid ideal.

More precisely,
\begin{theorem}\label{thm:main1}
The pullback  vector bundle $p^!A$ comes equipped with a “natural” almost Lie algebroid  inscribed in an exact sequence 
 $$   \xymatrix{ K\,\ar[d] \ar@{^(->}[r]& p^!A\ar@{->>}[r]\ar[d]
& \mathcal{D}_\mathcal{F}\ar[d]\\\mathrm{Nbl}_\mathcal{F}(M)\ar@{=}[r]&
\mathrm{Nbl}_{\mathcal{F}}(M)\ar@{=}[r]& \mathrm{Nbl}_{\mathcal{F}}(M)}$$
{where $\mathcal{D}_{\mathcal F}$ is a Lie algebroid such that $\mathcal{D}_{\mathcal F}|_{p^{-1}(M_{\mathrm{reg}})}\simeq \mathcal{F}|_{M_{\mathrm{reg}}}$, and $K$  is a vector bundle whose fiber at $V\in \mathrm{Nbl}_\mathcal{F}(M)$ is a subspace of $(p^!A)_V=A_{p(V)}$ stable under the pointwise skew-symmetric bracket of $A_{p(V)}$ and coincides with $\ker\rho_{p(V)}$ on regular points. By “natural” we mean that the projection $p^!A\to A$ is an almost Lie algebroid morphism. Also, $\mathcal{D}_\mathcal{F}\to\mathrm{Nbl}_{\mathcal{F}}(M)$ depends only on $\mathcal{F}$, not on $A$.}
\end{theorem}

While establishing the proof of Theorem \ref{thm:main1}, we observe the following outcome concerning Lie algebroids and singular foliation, for the case where $E=TM$.
 
The next statement gives a precise meaning to the slogan “A Lie algebroid splits in a purely-manifold part and a purely Lie algebra part”.  Recall that a Debord Lie algebroid can be thought as a “purely-manifold” Lie algebroid, in the sense that it is entirely described by its image through
the anchor map.

\begin{theorem}\label{thm:main}
Let $E=TM$. For any Lie algebroid $A\to M$, the pullback  vector bundle $p^!A$ comes equipped with a “natural” Lie algebroid structure inscribed in an exact sequence 
 $$   \xymatrix{ K\ar[d] \ar@{^(->}[r]& p^!A\ar@{->>}[r]\ar[d]
& \mathcal{D}_\mathcal{F}\ar[d]\\\mathrm{Nbl}_\mathcal{F}(M)\ar@{=}[r]&
\mathrm{Nbl}_{\mathcal{F}}(M)\ar@{=}[r]& \mathrm{Nbl}_{\mathcal{F}}(M)}$$
where $\mathcal{D}_{\mathcal F}$ is a Debord Lie algebroid such that $\mathcal{D}_{\mathcal F}|_{p^{-1}(M_{\mathrm{reg}})}\simeq \mathcal{F}|_{M_{\mathrm{reg}}}$, and $K$  is a Lie algebra bundle. By “natural” we mean that the projection $p^!A\to A$ is a Lie algebroid morphism. Also, $\mathcal{D}_\mathcal{F}\to\mathrm{Nbl}_{\mathcal{F}}(M)$ depends only on $\mathcal{F}$, not on $A$.
\end{theorem}

This makes sense of the following definition.
\begin{definition}
     The \emph{Nash blow-up} $\mathrm{Nash}(A)\to \mathrm{Nbl}_\mathcal{F}(M)$ of  an $E$-valued almost Lie algebroid/Lie algebroid $A\to M$ is the almost Lie algebroid/Lie algebroid structure on $p^!A$ of Theorem \ref{thm:main1}. \begin{enumerate}
         \item $\mathcal{D}_{\mathcal F}$ is the \emph{Lie algebroid/Debord Lie algebroid associated} to $\mathcal{F}=\rho(\Gamma(A))\subseteq \Gamma(E)$;
         \item $K$ is the \emph{tautological vector bundle/Lie algebra bundle associated} to $A$.
     \end{enumerate}
\end{definition}

The proof of Theorem \ref{thm:main1} and \ref{thm:main}  go through the assertions below.

Let us construct $K$ and $\mathcal{D}_\mathcal
F$. The short exact sequence of vector bundles \eqref{eq:shortsequence-vb1} on $\mathrm{Grass}_{-r}(A)$ restricts  to a short exact sequence of vector bundles, 
    \begin{equation}\label{eq:shortsequence}
    \xymatrix{
      K\ar@{^(->}[r] \ar[d]&\ar[d] p^! A \ar@{->>}[r] &  \mathcal{D}_{\mathcal F} \ar[d] \\\mathrm{Nbl}_\mathcal{F}(M)\ar@{=}[r]& \mathrm{Nbl}_\mathcal{F}(M)\ar@{=}[r] & \mathrm{Nbl}_\mathcal{F}(M)
 }
    \end{equation}

where  $K$ and by $\mathcal{D}_{\mathcal F}$  stand for the restriction of the  tautological bundle $\tau^{-r}_A$ (of codimension $r$) and  quotient bundle $Q_A^r$ (of rank $r$)  to $\mathrm{Nbl}_\mathcal{F}(M)$, respectively. 
\begin{proposition}\label{prop:main}In the notations above:
     
     \begin{enumerate}
     \item For every $a\in \Gamma(A)$, there exists a unique  vector field\,  $a^\dagger$ on $\mathrm{Nbl}_{\mathcal {F}}(M)$ such that $p_*(a^\dagger)=\rho(a)\circ p$. 
         \item The mapping $a\in \Gamma(A)\mapsto  a^\dagger\in\mathfrak X(\mathrm{Nbl}_\mathcal{F}(M))$ is a Lie algebroid action of $A$ on $\mathrm{Nbl}_\mathcal{F}(M)$, i.e., it satisfies for all $a,b\in \Gamma(A)$ and $f\in \mathcal{O}_M$
\begin{itemize}
     \item [i.] $[a,b]^\dagger=[a^\dagger,b^\dagger]$,
     \item [ii.] ${(fa)^\dagger}=(f\circ p)\,a^\dagger$,
     \item [iii.] $p_*(a^\dagger)=\rho(a)\circ p$.
 \end{itemize}
 \item  there is a unique  Lie algebroid structure on $p^!A$ such that the projection \begin{equation}
     \xymatrix{p^!A\ar[d]\ar[r]& A\ar[d]\\\mathrm{Nbl}_\mathcal{F}(M)\ar[r]_<<<<<p& M}
 \end{equation}
 is a Lie algebroid morphism and whose anchor map is given by $$\Hat{\rho}(p^!a)=a^\dagger,\; \;\text{for all $a\in \Gamma(A)$}.$$
 
 \item  $K\rightarrow \mathrm{Nbl}_\mathcal{F}(M)$ is a Lie algebra bundle and $K\subset p^!A$ is a Lie algebroid ideal whose fibers coincide with the fibers of  $\ker\rho$ on $p^{-1}(M_{\mathrm{reg}})$.
     \end{enumerate}
\end{proposition}
The following statement justifies the notation $\mathcal{D}_\mathcal{F}$. %Proposition \ref{prop:main} implies the following results.
\begin{corollary}\label{cor:main}
In the notations of Proposition \ref{prop:main}, 
there is a Lie algebroid structure $\left(\mathcal{D}_{\mathcal F},\,[\cdot, \cdot\,]_{\mathcal{D}_{\mathcal F}},\,\overline{\rho}\right)$ on the quotient bundle $\mathcal{D}_{\mathcal F}$. If $E=TM$, its anchor map $\overline{\rho}\colon \mathcal{D}_\mathcal{F}\to T\left(\mathrm{Nbl}_\mathcal{F}(M)\right)$ is injective on the open dense subset $p^{-1}(M_{\mathrm{reg}})$, i.e., it is a Debord Lie algebroid.
        %\item whose basic singular foliation $p^!\mathcal{F}$ is generated by the $a^\dagger$'s such that the following diagram commutes\begin{equation}\xymatrix{\mathcal{D}_{\mathcal F}\ar[r]^<<<<{\overline \rho} \ar[d]^{}& \ar[d]^{Tp}   T\left(\mathrm{Nbl}_\mathcal{F}(M)\right) \\A \ar[r]^{\rho} & TM    } \end{equation}\end{enumerate}
    
\end{corollary}

Proposition \ref{prop:main} and Corollary \ref{cor:main} implies that  \begin{equation}\label{}
        0\rightarrow K\rightarrow p^!A\rightarrow \mathcal{D}_{\mathcal F}\rightarrow 0
    \end{equation} is a short exact sequence of Lie algebroids. This proves Theorem \ref{thm:main}.

\begin{remark}
    For $E=TM$, the blow-up map $p\colon \mathrm{Nbl}_{\mathcal {F}}(M)\longrightarrow M$ sends the leaves of $\mathcal{D}_{\mathcal F}$ to the leaves of $\mathcal{F}$ and the restriction $p|_{p^{-1}(M_\mathrm{reg})}$ is an isomorphism of regular foliations. In particular, the Nash blow-up does not change the dimension of regular leaves.
\end{remark}

We end this section with a description of the leaves of $\mathcal{D}_\mathcal{F}$ with $E=TM$. Let $L\subseteq M$ be a leaf of the Lie algebroid $A\to M$. The restriction $A|_L\to TL$ is a transitive Lie algebroid on  $L$. The latter acts on $p^{-1}(L)\hookrightarrow \mathrm{Grass}_{-(r-\dim L)}(\ker \rho|_L)$. The orbits of this action are those of $\mathcal{D}_\mathcal{F}|_{p^{-1}(L)}$ by construction. Let us describe these leaves.\\

For simplicity, we assume that $A$ is integrable to a source-simply connected Lie groupoid $\Gamma\rightrightarrows M$. For $L\subseteq M$, we shall denote $\Gamma^L_L:=s^{-1}(L)\cap t^{-1}(L)$ and $\Gamma^m_m:=\Gamma^{\{m\}}_{\{m\}}$, where $s,t$ are source and target of $\Gamma$. 
\begin{proposition}\label{prop:leaves-description}
     Let $p^!\mathcal{F}$ be the basic singular foliation of $\mathcal{D}_\mathcal{F}$ and $L\subset M$ a leaf of $\mathcal{F}$. The restriction groupoid $\Gamma^L_L\rightrightarrows L$ to $L$  %The simply connected Lie group $G_m(\mathcal{F})$  integrating $\mathfrak g_m(\mathcal{F})$ 
     naturally acts on $\mathrm{Grass}_{-(r-\dim L)}(\ker \rho|_L)$. 

     \begin{enumerate}
         \item this action preserves the fiber $$p^{-1}(L)\hookrightarrow\mathrm{Grass}_{-(r-\dim L)}(\ker \rho|_L)$$ of the blow-up map $p\colon \mathrm{Nbl}_\mathcal{F}(M)\rightarrow M$. 
         \item If $L=\{m\}$ is a leaf of $\mathcal{F}$, then the leaf $L_{V_m}$ of $p^{!}\mathcal{F}$ that passes through $V_m\in p^{-1}(m)$ corresponds to a $\Gamma^L_L$-orbit and we have $L_{V_m}\simeq \frac{\Gamma^m_m}{\mathrm{Stab}(V_m)}$. Moreover,
\begin{enumerate}
    \item $V_m\in p^{-1}(m)$ is a leaf of $p^!\mathcal{F}$ if only if $V_m$ is a Lie ideal of $\mathfrak g_m(\mathcal{F})$.
    \item $V_m\in p^{-1}(m)$ is a regular point of $p^!\mathcal{F}$ if only if $\mathrm{Stab} (V_m)=V_m$.
\end{enumerate}
where $\mathrm{Stab}(V_m)$ is the stabilizer of $V_m\in p^{-1}(m)$. In that case, the Debord Lie algebroid $\mathcal{D}_\mathcal F|_{p^{-1}(m)}$ is $$\dfrac{\Gamma^m_m\times  \mathfrak g_m(\mathcal F)/{V_m}}{N(V_m)}\to \frac{\Gamma^m_m}{N(V_m)}$$  where $N(V_m)$ is the normalizator of $V_m\in p^{-1}(m)$.
     \end{enumerate}
\end{proposition}

\subsubsection{Proofs of main results}\label{sec:main-statement}
In this section, we prove Proposition \ref{prop:main} and Corollary \ref{cor:main}. The latter are proved through Proposition \ref{prop:lift:1}, Lemma \ref{grass_pullback} and Lemma \ref{lem:anchorGrass}.

Here,  we write the proofs in the context $E=TM$, that is, for almost Lie algebroids over singular foliations $\mathcal{F}\subseteq \mathfrak X(M)$. The singular subalgebroid case is similar.

\begin{lemma}[Lift of vector fields of $\mathcal{F}$ on $\mathrm{Nbl}_\mathcal{F}(M)$] 
\label{grass_pullback}In the notations of Proposition \ref{prop:lift:1},
  consider a lift $ \hat{\xi}_X$ of some $ X \in \mathcal F$ to the Grassmann bundle $\mathrm{Grass}_{-r}(A)$ of some anchored bundle $(A,\rho) $, where $\xi_X$ is a linear vector field on $A$ over $X$. 
\begin{enumerate}
\item The flow $\phi_t^{\hat{\xi}_X}$ of $ \hat{\xi}_X$, whenever it is defined, preserves the subset $\mathrm{Nbl}_\mathcal{F}(M)\subset \mathrm{Grass}_{-r}(A)$, i.e. 
    \begin{equation}
    \xymatrix{\mathrm{Nbl}_\mathcal{F}(M)\ar[r]^{\phi_t^{\hat{\xi}_X}}\ar[d]_{p}&\mathrm{Nbl}_\mathcal{F}(M)\ar[d]^{p}\\M\ar[r]^{\phi_t^X}&M}
\end{equation}
\item In particular, $\hat{\xi}_X$ is tangent to $\mathrm{Nbl}_\mathcal{F}(M)$ is the neighborhood of any point where $\mathrm{Nbl}_\mathcal{F}(M)$  is a sub-manifold. 
\end{enumerate}
 \end{lemma}

 \begin{proof}
    For every $x\in M$, the flow $\phi_t^{\xi_X}$ (whenever it makes sense) sends $\mathrm{\ker}{\rho}_{x}$ to  $\ker{\rho}_{\phi^X_t(x)}$, that is
$$\phi_t^{\xi_X}|_{x}\left(\mathrm{\ker}{\rho}_{x}\right)=\ker{\rho}_{\phi^X_t(x)}$$in view of Diagram \eqref{eq:diag1}. In particular, in terms of the section $\sigma_{\ker}: M_{\mathrm{reg}} \to \mathrm{Grass}_{-r}(A)$,\, $x\mapsto \ker\rho_x$, it means that
 $$ \phi^{\hat{\xi}_X}_t\circ \sigma_{\ker}=\sigma_{\ker} \circ \phi^X_t.$$
 In particular, since $\phi^X_t(M_{\mathrm{reg}})\subseteq M_{\mathrm{reg}}$, this implies that $  \phi^{\hat{\xi}_X}_t $ preserves the closure $ \overline{\sigma_{\ker}(M_{\mathrm{reg}})}$, i.e., it preserves $\mathrm{Nbl}_\mathcal{F}(M)$. This completes the proof.
 \end{proof}

\begin{lemma}
\label{lem:anchorGrass}
In the smooth case, there exists a vector bundle morphism $$\xymatrix{\Tilde{\rho}\colon P^!A\ar[r]\ar[d]& T(\mathrm{Grass}_{-r}(A))\ar[d]\\\mathrm{Grass}_{-r}(A)\ar@{=}[r]&\mathrm{Grass}_{-r}(A)
}$$ 
such that
 \begin{enumerate}
 \item The following diagram commutes:
$$ \xymatrix{P^!A \ar[r]^<<<<{\Tilde\rho} \ar[d]^{}& \ar[d]^{TP}   T \left(\mathrm{Grass}_{-r}(A)\right) \\A \ar[r]^{\rho} & TM    } $$
  \item For any $a \in \Gamma(A)$, $\Tilde{\rho}(P^!a) \in \mathfrak X\left(\mathrm{Grass}_{-r}(A)\right)$ is a lift of $X=\rho(a)\in \mathcal{F}$ tangent to $\mathrm{Nbl}_{\mathcal{F}}(M)$. Here $P^!a $ is the pullback  of the section  $a$. 
 \end{enumerate} 
 Here, $P^!A$ stands for the pullback  of $A$ on the manifold $\mathrm{Grass}_{-r}(A))$ and $P^!A\rightarrow A$ is the projection onto $A$. {In the holomorphic/real analytic cases, existence is local.}
\end{lemma}
\begin{proof}
Let $V_x \in{\mathrm{Grass}}_{-r}(A)\to M$ be an element of the fiber at $x\in M$. Consider a local trivialization $a_1, \dots, a_d$ of $A$ in a neighborhood $\mathcal U $ of $x$. {The evaluation at $V_x\in {\mathrm{Grass}}_{-r}(A)$ of the lifts $\hat{\xi}_{a_1}, \dots,\hat{\xi}_{a_i}:= \hat{\xi}_{{\rho}(a_i)},\cdots,\hat{\xi}_{a_d}$ constructed as in Lemma \ref{grass_pullback} defines a vector bundle morphism $\Tilde{\rho}_{\mathcal U} $ that satisfies item 1 and item 2, by construction. In the smooth case, using a partition of unity $ (\mathcal{U}_i,\phi_i,\Tilde{\rho}_{\mathcal{U}_i})_{i \in I}$, this construction can be made global: it is easily checked that  $ \sum_i P^*\phi_i  \Tilde{\rho}_{\mathcal{U}_i} $  still satisfies item 1 and item 2, because each $\Tilde{\rho}_{\mathcal{U}_i}$ does.}
\end{proof}

\noindent
\textbf{Proof of Proposition \ref{prop:main} and  \ref{prop:leaves-description} }

\begin{proof}[Proof of Proposition \ref{prop:main}]
\begin{enumerate}
    \item Let $a\in \Gamma(A)$ and consider $\rho(a)=X$. In the notations of Lemma \ref{grass_pullback}, let $\hat{\xi}_X $   be the lift of $X$ to $\mathrm{Grass}_{-r}(A)$. The restriction of $\hat{\xi}_X $ to ${\mathrm{Nbl}}_\mathcal F(M) $ is the required vector field. Uniqueness is a direct consequence of the fact that $p$ is an isomorphism when restricted to the dense open subset $p^{-1}(M_{\mathrm{reg}})$.

    \item The construction of the map  $a\mapsto a^\dagger$ is obtained by restricting the vector bundle morphism $$P^!A\stackrel{\title{\rho}}{\rightarrow} T\left(\mathrm{Grass}_{-r}(A)\right)$$ of Lemma \ref{lem:anchorGrass} to $\mathrm{Nbl}_\mathcal{F}(M)$. The vector fields $[a, b]_A^\dagger$ and $[a^\dagger, b^\dagger]$ are both $p$-related to $[\rho(a), \rho(b)]$, for all $a, b\in \Gamma(A)$. By uniqueness, we have
 $$ [a, b]_A^\dagger=[a^\dagger, b^\dagger], \, \, \hbox{for all $a, b\in \Gamma(A)$}.$$ 

\item  The Lie algebroid structure on $p^!A$ is obtained as follows: the anchor map $\Hat\rho$ is defined by $p^!a\mapsto a^\dagger$ and  the Lie bracket on $p^!A$ as $$[p^!a, p^!b]_{p!A}:= p^!([a,b]_A)$$ for all $a, b\in \Gamma(A)$. We extend by Leibniz identity to all sections.
\item The fact that $K\rightarrow \mathrm{Nbl}_\mathcal{F}(M)$ is a Lie algebra bundle is a direct consequence of item 1 of Proposition \ref{prop:almost}. Let us show that it is a Lie algebroid ideal of $p^!A$.  Notice that  $\Gamma(K)\subset \ker \Hat{\rho}$ and $\Gamma(K)|_{p^{-1}(M_\mathrm{reg})}=\ker \Hat{\rho}|_{p^{-1}(M_\mathrm{reg})}$ %the fibers of  $K$ are included in the kernels of the anchor map $\Hat{\rho}$ of $p^!A$
: indeed, we have $\ker\rho_x\subseteq \ker\Hat{\rho}_{V_x}$ for all $V_x\in p^{-1}(M_{\mathrm{reg}})$ by construction, and both spaces coincide on $p^{-1}(M_{\mathrm{reg}})$ since $\dim (\mathrm{Im}\Tilde{\rho}_{V_x})=\dim(\mathrm{Im}{\rho}_{x})$ for all $V_x\in p^{-1}(M_{\mathrm{reg}})$. For $E\neq TM$, we only have the inclusion $\Gamma(K)\subset \ker \Hat{\rho}$.

{Now, take $\xi\in \Gamma(K)$, $\ell\in \Gamma(p^!A)$ and $V\in \mathrm{Nbl}_\mathcal F(M)$. Let $p^!\rho\colon p^!A \to p^!TM$ be the pullback of the anchor map $\rho\colon A\to TM $ along $p\colon \mathrm{Nbl}_\mathcal{F}(M)\to M$. One has, \[p^!{\rho}\left([\xi, \ell]_{p^!A}\right)=0, \;\;\text{on $p^{-1}({M}_{\mathrm{reg}})$}.\]
Let $V_n\in p^{-1}(M_{\mathrm{reg}})$ be a sequence that converges to $V$. For all $n\in \mathbb N$, we have $$[\xi, \ell]_{p^!A}(V_n)\in \ker (p^!{\rho})_{V_n}=\ker{\rho}_{p(V_n)}=V_n.$$ This implies that the limit $[\xi, \ell]_{p^!A}(V_n)\underset{n \to +\infty}{\longrightarrow} [\xi, \ell]_{p^!A}(V) $ belongs to $V$. This proves that $$\left[\Gamma(K),\; \Gamma(p^!A)\right]_{p^!A}\subseteq \Gamma(K).$$
This completes the proof.}
 \end{enumerate}
\end{proof}
 \begin{proof}[Proof of Corollary \ref{cor:main}]
We use the notations of Proposition  \ref{prop:main}.
By Proposition \ref{prop:main}(4), $K\subset p^!A$ is a Lie algebroid ideal, therefore the anchor map $\Hat{\rho}\colon p^!A\mapsto T \left(\mathrm{Nbl}_\mathcal{F}(M)\right)$ goes to the quotient 
\begin{equation}
 \xymatrix{0\ar[r] &K\ar[r]&p^!A\ar[r]\ar[d]_{\Hat{\rho}} &\mathcal{D}_{\mathcal F}=p^!A/K\ar[r]\ar@{-->}[ld]^{\overline{\rho}}&0\\& &T(\mathrm{Nbl}_\mathcal{F}(M)) & &.}
 \end{equation} to equip $\mathcal{D}_{\mathcal F}$ with a Lie algebroid structure. If $E=TM$, the fibers of the tautological bundle $K\to \mathrm{Nbl}_\mathcal{F}(M)$ coincide with the fibers $\{\ker\Hat{\rho}_{V_x}\}_{V_x\in p^{-1}(x), x\in M}$ on the open dense subset  $p^{-1}(M_{\mathrm{reg}})$, thus,  $\mathcal{D}_{\mathcal F}$ is a Debord Lie algebroid.%injective on the open dense subset  $p^{-1}(M_{\mathrm{reg}}) $ with image $p^!\mathcal F $.    Therefore, the singular foliation $p^!\mathcal{F}$ is  projective  and is the image of the induced Lie algebroid structure  on  $\mathcal{D}_{\mathcal F}$. 
\end{proof}

\begin{proof}[Proof (of Proposition \ref{prop:leaves-description})]
    We give a proof with a leaf $L=\{m\}$ reduced to a point, the general case is similar. Let $(A, [\cdot\,,\cdot]_A, \rho)$ be a Lie algebroid and $m\in M$. Notice that for every $a\in \Gamma(A)$, the flow $\phi^{\xi_X}_t\colon A\rightarrow A$ of the linear vector field $\xi_X\in \mathfrak X(A)$ associated to $\mathcal{L}_a$, with $\rho(a)=X$, satisfies for all $b\in \Gamma(A)$ \begin{equation}\label{eq:formula}
        \frac{d}{dt}_{|_{t=0}}(\phi^{\xi_X}_t)^*b=[a, b]_A
    \end{equation}with $$(\phi^{\xi_X}_t)^*(b)(m):=\phi^{\xi_X}_{-t}(m)(b(\phi^X_t(m))),$$ whenever the flows make sense. Using the formula \eqref{eq:formula}, we make the following observation: let $G_m(\mathcal{F})$ be the simply connected Lie group integrating $\mathfrak g_m(\mathcal{F})$. The latter naturally acts on 
$\mathrm{Grass}_{-(r-\dim(L_m)}\left(\mathfrak g_m(\mathcal F)\right)$ through the adjoint action. By Lemma \ref{grass_pullback} (1),  this action preserves the fiber $$p^{-1}(m)\hookrightarrow\mathrm{Grass}_{-(r-\dim(L_m)}\left(\mathfrak g_m(\mathcal F)\right)$$ of the blow-up map $p\colon \mathrm{Nbl}_\mathcal{F}(M)\rightarrow M$. If $m\in M$ is a leaf of $\mathcal{F}$, then the leaf $L_{V_m}$ of $p^{!}\mathcal{F}$ that passes through $V_m\in p^{-1}(m)$ is a ${G}_m(\mathcal{F})$-orbit, since the flows $(\phi^{\xi_X}_t)_{X\in \mathcal{F}}$ generate the foliation $p^!\mathcal{F}$. It is easily checked that
\begin{enumerate}
    \item $V_m\in p^{-1}(m)$ is a leaf of $p^!\mathcal{F}$ if only if $[V_m]$ is a Lie ideal of $\mathfrak g_m(\mathcal{F})$.
    \item Since $$\dim L_{V_m}=r- \frac{\mathrm{Stab} ([V_m])}{[V_m]},$$  $V_m\in p^{-1}(m)$ is a regular point if only if $\mathrm{Stab} ([V_m])=[V_m]$.
\end{enumerate}
\end{proof}

\section{Examples}\label{sec:3}
In this section, we present some examples related to singular foliations, Lie algebroids, and also to Poisson manifolds. The main results of the paper in Section \ref{sec:2} provide a desingularization method for $E$-valued (almost) Lie algebroids over a manifold $M$. They can also be used to desingularize finitely generated singular subalgebroids $\mathcal{F}\subseteq \Gamma(E)$, particularly finitely generated singular foliations $\mathcal{F}\subseteq\mathfrak X(M)$ by choosing an almost Lie algebroid $A$ over $\mathcal{F}$. The blow-up space $\mathrm{Nbl}_\mathcal{F}(M)$ and the Debord foliation $\mathcal{D}
_\mathcal{F}$ as in Theorem \ref{thm:main1} are independent of the choice of the almost Lie algebroid by Proposition \ref{prop:well-defined}. Our primary focus lies in computing the blow-up space $\mathrm{Nbl}_\mathcal{F}(M)$ and the Debord Lie algebroid $p\colon \mathcal{D}
_\mathcal{F}\to M$ associated with $\mathcal{F}$. We shall denote the basic singular foliation of $\mathcal{D}
_\mathcal{F}$ by $p^!\mathcal{F}$.

\subsection{Lie algebroids, Lie algebra actions, singular foliations}
\begin{example}[The Nash blow-up does not affect projectivity]
    %If $(M, \mathcal{F})$ is a Debord singular foliation of rank $r$,  i.e., there is a 
    If $(A, [\cdot\,,\cdot]_A,\rho)$ is a Lie algebroid whose anchor map is injective on an open dense subset, then $$ \xymatrix{ \mathcal{D}_\mathcal{F}\ar[r]^{\simeq}\ar[d]& A\ar[d]\\\mathrm{Nbl}_\mathcal{F}(M)\ar[r]^<<<<{\simeq}&M}$$ since $\mathrm{Grass}_{-\mathrm{rk}(A)}(A)\simeq M$.

    For instance, the Nash blow-up has no effect on the Lie algebroid of  vector fields on a $b$-manifold in the sense of \cite{GUILLEMIN201486}. 
\end{example}

%\begin{example}[Singular subalgebroids]Let $(E, [\cdot\,,\cdot]_E, \rho_E)$ be a Lie algebroid whose regular leaves are of dimension $r$ and   $\mathcal{F}\subseteq \Gamma(E)$  a singular subalgebroid as in Definition \ref{def:Sing-subal}. Assume that that blow-up space $p\colon\mathrm{Nbl}_{\mathcal{F}_E}(M)\rightarrow M$ is smooth with $\mathcal{F}_E:=\rho_E(\Gamma(E))$. In the notations of Proposition \ref{prop:main}, the pullback $p^!\mathcal{F}$ of $\mathcal{F}$ on $\mathrm{Nbl}_{\mathcal{F}_E}(M)$ is a singular subalgebroid of $p^!E$. If $\ker\rho_E\subseteq \mathcal F$ then $\Gamma(\overline{\tau}^r_E)\subset p^!\mathcal{F}$ and $p^!\mathcal{F}/\Gamma(\overline{\tau}^r_E)$ is a Lie-Rinehart subalgebra of $\Gamma(\overline{Q}_E^r)$ whose anchor map in injective. ince the anchor map of the  Lie algebroid structure on $\overline{Q}_E^r$ is almost injective, it is integrable to a Lie groupoid over $\mathrm{Nbl}_{\mathcal{F}_E}(M)$ by \cite{Debord}. Hence, the notion of bisubmersion and holonomy  groupoid are well-defined for $p^!\mathcal{F}/\Gamma(\overline{\tau}^r_E)$, see \cite{MZ}.\end{example}

\begin{example}[Vector fields tangent to a closed submanifold]\label{ex:CDO} Let $M$ be a manifold and $W \subseteq M$ a closed submanifold. We denote by $\mathcal{F}_W=\mathfrak X_W (M)$ the singular foliation of all vector fields tangent to $W$ \cite{LLL1}. It is still an open question whether $\mathcal{F}_W$ is the image through the anchor map of a Lie algebroid on the whole manifold $M$. However, the answer is positive in a tubular neighborhood of $W$ in $M$. Indeed, such neighborhood comes equipped with a diffeomorphism to a neighborhood of the zero section of the normal bundle $\mathcal{N}_W= \frac{TM|_W}{TW}$. Under this diffeomorphism,
the restriction of $\mathcal{F}_W$ to that neighborhood corresponds to the space of all vector fields on $\mathcal{N}_W$ tangent to the zero section, which we denote again by $\mathcal{F}_W$. It follows that:
\begin{enumerate}
    \item Linear vector fields on $\mathcal{N}_W$ are the sections of a transitive Lie algebroid on  $W$  known as $\mathrm{CDO}(\mathcal{N}_W)$ in \cite{Mackenzie},
\item  $\mathrm{CDO}(\mathcal{N}_W)$ acts on $\mathcal{N}_W$,
\item  $\mathcal F_W$ is the image the anchor map of this action Lie algebroid on $\mathcal{N}_W$.
\end{enumerate}
For simplicity, we can assume that the tubular neighborhood is $M$. It is easy to check that $p\colon \mathrm{Nbl}_{\mathcal{F}_W}(M)\rightarrow M$ is the usual blow-up of $M$ along $W$ (see, e.g. \cite{Huybrechts}) and $p^{-1}(W)\simeq \mathbb P(\mathcal{N}_W)$. The pullback $p^!\mathcal{F}_W\simeq \mathcal{D}_{\mathcal{F}_W}$ of $\mathcal{F}$ along $p\colon \mathrm{Nbl}_{\mathcal{F}_W}(M)\rightarrow M$ is the Debord singular foliation made of all vector fields on $\mathrm{Nbl}_{\mathcal{F}_W}(M)$ that are tangent to $p^{-1}(W)$.

In conclusion, we have the short exact sequence of Lie algebroids over the blow-up of $\mathbb K^d$ at $0$ as follows $$0\to K\to p^!(\mathrm{CDO}(\mathcal{N}_W))\to \mathcal{D}_{\mathcal{F}_W}\to 0, $$
here $K$ is made of pairs $(\ell,\mathfrak M)\in \mathbb P(\mathcal{N}_W)\times \mathfrak {gl}(\mathcal{N}_W)$ such that $\ell$ is an eigenspace of $\mathfrak M$.
\end{example}

The following example is a particular case of Example \ref{ex:CDO} with $W$ being a single point.
\begin{example}[Action of $\mathfrak {gl}_d(\mathbb K)$ on $M=\mathbb K^d$]\label{ex:gln_1} Consider the action Lie algebroid $A=\mathfrak {gl}_d(\mathbb K)\ltimes \mathbb K^d\to \mathbb K^d$. Its basic singular foliation  $ \mathcal F=\mathrm{Span}(X_{kj}=x_k\frac{\partial}{\partial x_j},\; k,j=1,\ldots,d)$ is generated by all vector fields vanishing at the origin. It is easily checked that the Nash blow-up $p\colon \mathrm{Nbl}_\mathcal{F}(M)\to \mathbb K^d$ coincides to the usual blow-up of $\mathbb K^d$  at the origin. In the $i$-th coordinate chart $(\mathcal{U}_i, y_1\ldots,y_d)$ of the blow-up space, the pullback $p^!\mathcal{F}$ is generated by the $\frac{\partial}{\partial y_j}$'s with $j\neq i$ and $y_i\frac{\partial}{\partial y_i}.$ These computations show that $p^!\mathcal{F}\simeq \mathcal{D}_\mathcal{F}$ is generated by all vector fields tangent to the projective space $p^{-1}(0)=\mathbb P^{d-1}$.

We have the short exact sequence of Lie algebroids $$0\to K\to p^!(\mathfrak {gl}_d(\mathbb K)\ltimes \mathbb K^d)\to \mathcal{D}_\mathcal{F}\to 0, $$
here $K$ is made of pairs $(\ell,\mathfrak M)\in \mathbb P^{d-1}\times \mathfrak {gl}_d(\mathbb K)$ such that $\ell$ is an eigenspace of $\mathfrak M$.
\end{example}

\begin{example}[Vector fields vanishing at order $k\in \mathbb N$ on $M=\mathbb K^d$]
Denote by $x_1, \ldots, x_d$ the canonical coordinates of $\mathbb K^d$. Let $k\in \mathbb N$ and  $\mathcal{F}^k$ be the singular foliation on $\mathbb K^d$ generated by the vector fields $x^I\frac{\partial}{\partial x_j}$ for $j=1, \ldots, d$ and $x^I$ ranges over all monomials of degree $k$ in $x_1, \ldots, x_d$. For every $k\in \mathbb N$, $M_{\mathrm{reg}}= \mathbb K^d\setminus \{0\}$. %$x^2\frac{\partial}{\partial x},\, xy\frac{\partial}{\partial x}, y^2\frac{\partial}{\partial x},\, x^2\frac{\partial}{\partial y},\, xy\frac{\partial}{\partial y}, $  and $y^2\frac{\partial}{\partial y}$. 
The singular foliation $\mathcal{F}^k$ is the image of an almost Lie algebroid $A=\mathbb K^d\times \mathbb K^{N}$ whose anchor map $\rho$ is $e_{I,j}\mapsto x^I\frac{\partial}{\partial x_j}$ and the almost Lie bracket is given by \begin{equation}
    [e_{I,j}, e_{J,l}]_A:=x^{J\setminus j}e_{I,l}- x^{I\setminus l}e_{J,j}
\end{equation}
with $N=d\times\begin{pmatrix} d+k-1\\k\end{pmatrix}$ and $x^{I\setminus j}$ is the monomial $\frac{x^I}{x_j}$. For $k=1$, $\mathcal{F}^1$ is given by the action $\mathfrak{gl}_d(\mathbb K) \curvearrowright\mathbb K^d$, see Example \ref{ex:gln_1}. For $k\neq 1$, it is an open question whether $\mathcal{F}^k$ is the image of some Lie algebroid or not,  see \cite{LLS}.

A direct computation shows that for every $k\in \mathbb N$,  $$\mathrm{Nbl}_{\mathcal{F}^k}(M)\simeq\mathrm{Bl}_{\mathcal{I}_k}(\mathbb K ^d).$$where $\mathrm{Bl}_{\mathcal{I}_k}(\mathbb K^d)$ is the usual blow-up of $\mathbb K^d$ along the ideal $\mathcal{ I}_k$ of all polynomial functions vanishing at order $k$. Consider the embedding \begin{equation}
    \label{eq:veronese}\nu_k\colon \mathbb P^{d-1}\rightarrow \mathbb P^{N-1},\, [x_1:\cdots:x_d]\mapsto [\cdots: x^I: \cdots ]
\end{equation} known as the Veronese map of degree $k$, see e.g., \cite{harris}. For $k=2$ and $d=2$, the map \eqref{eq:veronese} reads $$\nu_2\colon \mathbb P^{1}\rightarrow \mathbb P^{2},\, [x:y]\mapsto [x^2:xy :y^2 ].$$ %For small integers, $\nu_k(\mathbb P^{1})$ is known as the rational normal curve; $\nu_2(\mathbb P^{1})$ is the standard parabola in coordinates $(t, t^2)$; $\nu_3(\mathbb P^{1})$ is the twisted cubic in coordinates $(t, t^2, t^3)$ see e.g \cite{harris}. 
Under the embedding \eqref{eq:veronese}, $\mathrm{Bl}_{\mathcal{I}_k}(\mathbb K^d)$ is isomorphic to the usual blow-up of $\mathbb K^d$ at $0$. Also, $p_k^{-1}(0)$ is isomorphic to  $\nu_k(\mathbb P^{d-1})$, where $p_k\colon \mathrm{Nbl}_{\mathcal{F}^k}(M)\rightarrow M$ is the projection map. 
\end{example}

\begin{example}[$\mathfrak{sl}_2\curvearrowright \mathbb R^2$]
    Consider the Lie algebra $\mathfrak{sl}_2\subset \mathfrak {gl}_2$ generated by the canonical generators $h, e, f$ that satisfy $[h, e] =
2e$, $[h, f ] = -2f$ and $[e, f ] = h$. It acts on $\mathbb  R^2$ as: $\rho\colon  \mathbb R ^2\times\mathfrak{sl}_2\rightarrow T\mathbb R^2$

$$\rho\begin{pmatrix}
    h\\e\\f
\end{pmatrix}(x,y)=\begin{pmatrix}
    x&-y\\0&x\\y&0
\end{pmatrix}\begin{pmatrix}
    \frac{\partial}{\partial x}\\\frac{\partial}{\partial y}
\end{pmatrix}$$
Here, $x, y$ stand for the coordinates of $M=\mathbb R^2$. The set of regular points is  $M_{\mathrm{reg}}=\mathbb R^2\setminus \{0\}$. The blow-up space  $p\colon \mathrm{Nbl}_\mathcal{F}(M)\subset \mathbb R^2\times \mathrm{Grass}_{-2}(\mathfrak{sl}_2)\to \mathbb R^2$ is the usual blow-up of $\mathbb R^2$ along the ideal $I=(-x^2, xy, y^2)$. This corresponds to the usual blow-up $\mathrm{Bl}_0(\mathbb R^2)$ of $\mathbb R^2$ at $0$: indeed, we have $\nu\left(\mathrm{Bl}_0(\mathbb R^2)\right)=\mathrm{Nbl}_\mathcal{F}(M)$ under the closed embedding \begin{equation}
    \nu\colon \mathbb R^2\setminus\{0\}\times \mathbb P^1\hookrightarrow \mathbb R^2\setminus\{0\}\times \mathbb P^2,\;\left((x,y),[x:y]\right)\mapsto \left((x,y),[-x^2:xy: y^2]\right).
\end{equation}  The blow-up space $\mathrm{Nbl}_\mathcal{F}(M)$ is smooth. % by applying directly Proposition \ref{prop:smoothness}: the kernel $\ker \rho$ is a projective $\mathcal{O}_{\mathbb R^2}$-submodule: It is the line bundle generated by the section $xy h + y^2e -x^2f$, see Example 3.31 of \cite{LLS}. 
By a direct computation one can find the generators of the Debord singular foliation $p^!\mathcal{F}$, e.g., on the $x$-chart $\mathcal{U}_x$ of the blow-up space, it is generated by the vector fields $\overline h= x \frac{\partial}{\partial x}-2y\frac{\partial}{\partial y}$,\; $\overline e=\frac{\partial}{\partial y}$ and $\overline {f}=xy\frac{\partial}{\partial x}-y^2\frac{\partial}{\partial y}$. Since $\overline f=y\overline h+y^2\overline e$, the Debord singular foliation $p^!\mathcal{F}$ is generated by the $\mathcal{O}_{\mathbb R^2}$-linearly independent vector fields $\overline h, \,\overline e$ on $\mathcal{U}_x$. Here, $p^{-1}(0)=\{x=0\}$ is the unique singular leaf of $p^!\mathcal{F}$.
\end{example}
\subsection{Foliated bivector fields}\label{application:poisson}
In this section, we give examples of our blow-up construction in Poisson geometry.

For a bivector field $\pi=\{\cdot\,,\cdot\}\in \mathfrak X^2(M)$ we denote by $\pi^\#\colon T^*M\longrightarrow TM$ the vector bundle morphism which is defined by the relation $\langle\pi^\#(\alpha), \beta\rangle=\langle\pi, \alpha\wedge\beta\rangle$ with $\alpha, \beta \in \Omega^1(M)$.

The result of this section can be applied to a more general concept than Poisson structures. Let us recall the following definition introduced in \cite{TURKI201571}.
\begin{definition}
    A bivector field $\pi\in \mathfrak X^2(M)$ is said to be \emph{foliated} if the submodule $\pi^\#(\Omega^1(M))\subset \mathfrak X(M)$ is a singular foliation, i.e., \begin{equation}
        \left[\pi^\#(\Omega^1(M)), \pi^\#(\Omega^1(M))\right]\subset \pi^\#(\Omega^1(M)).
    \end{equation} 
    The singular foliation $\pi^\#(\Omega^1(M))$ shall be denoted by $\mathcal{F}_\pi$ and be called the \emph{basic singular foliation of $\pi$}.\end{definition}
\begin{example}Poisson structures on manifolds \cite{CPA,Crainic2021LecturesOP} are the first class of  examples of foliated bivector fields. Recall that a Poisson structure on a manifold $M$ is a bivector field $\pi=\{\cdot\,,\cdot\}\in \mathfrak X^2(M)$ fulfilling $[\pi, \pi]_{\mathrm{NS}}=0$ where $[\cdot\,,\cdot]_{\mathrm{NS}}$ is the Schouten–Nijenhuis bracket on multivector fields. For every function $h\in \mathcal O_M$ we associate a vector field $X_h:=\pi^\#(dh)=\{\cdot\,,h\}\in\mathfrak X(M)$ called the \emph{Hamiltonian} vector field of $h$. One has that $[X_h, X_g]=X_{\{g,h\}}$ for all $g,h\in \mathcal O_M$. This implies that every Poisson structure $\pi$ on a manifold is foliated. Another class of examples of foliated bivector fields are twisted-Poisson structures, i.e.,$[\pi,\pi]_{\mathrm{NS}}=\wedge^3\pi^\#(\alpha)$, where $\alpha\in\Omega^3(M)$  is a closed 3-form \cite{SW,TURKI201571}.
\end{example}

%For more details and examples of foliated bivector fields, see \cite{TURKI201571}.  

    For any foliated bivector field $(M, \pi)$, $\pi^\#\colon T^*M\rightarrow TM$ is an anchored bundle over $\mathcal{F}_\pi$ and thus can be equipped with an almost Lie bracket structure $[\cdot\,,\cdot]_\pi$. From now on, we consider $\pi^\#\colon T^*M\rightarrow TM$ with an almost Lie algebroid structure  over $\mathcal F_\pi$. When $\pi$ is Poisson \cite{Crainic2021LecturesOP} or twisted Poisson\cite{SW}, in both cases there exists Lie algebroid structure on cotangent bundle $T^*M$ whose anchor map is $\pi^\#\colon T^*M\rightarrow TM$. Therefore, Theorem \ref{thm:main} is applicable in all these cases.\\

In the notations of Section \ref{sec:Nbl-construction}, the following proposition shows that for every foliated bivector field $(M, \pi)$ the Nash construction $\mathrm{Nbl}^{\mathrm{Im}}_{\mathcal{F}_\pi}(M)\subset\mathrm{Grass}_{r}(TM) $ out of the image   and the one $\mathrm{Nbl}^{\ker}_{\mathcal{F}_\pi}(M)\subset\mathrm{Grass}_{-r}(T^*M)$  out of the kernel  are canonically isomorphic.

\begin{proposition}\label{prop:Poisson}
    Let $(M, \pi)$ be a (foliated) bivector field  and  $\pi^\#\colon T^*M\rightarrow TM$ its associated almost Lie algebroid over  $\mathcal{F}_\pi$. Then $\mathrm{Nbl}^{\ker}_{\mathcal{F}_\pi}(M)\simeq\mathrm{Nbl}^{\mathrm{Im}}_{\mathcal{F}_\pi}(M)$.
\end{proposition}

\begin{proof}
%The proof is a direct consequence of the following facts. It is well-known that 
There is a canonical bijection between $\mathrm{Grass}_{-(d-r)}(TM)$ and  $\mathrm{Grass}_{-r}(T^*M)$ which maps,  $V\subset T_xM$ to its annihilator $V^\circ\subset T^*_xM$, with $d=\dim M$.  On the other hand, since $\pi$ is skew-symmetric, i.e.,  $0=\langle\pi^\#(\alpha), \beta\rangle=-\langle\alpha, \pi^\#(\beta)\rangle$ for all $\alpha\in \ker(\pi^\#_x)$ and $\beta\in T^*_xM$,  the annihilator of $\mathrm{Im}(\pi^\#_x)$ is $\ker(\pi^\#_x)$
\end{proof}

%A corollary of Proposition \ref{prop:main} is the following,
%\begin{proposition}Let $(M, \pi)$ be a foliated bivector field and $\pi^{\#}\colon T^*M\rightarrow M$ be the associated almost Lie algebroid on $M$.  \begin{enumerate}\item there is a Lie algebroid $(\mathrm{NBl}(T^*M),[\cdot, \cdot\,]_{\mathrm{NBl}(T^*M)}, {\rho}^!)$ on $\mathrm{NBl}(M,\mathcal{F}_\pi)$ whose anchor map is injective on $p^{-1}(M_{\mathrm{reg}})$ and whose basic singular foliation $\mathrm{NBl}(\mathcal{F}_\pi)$ is generated by the $\mathrm{NBl}(X_f)$'s for $f\in \mathcal O_M$.\item In the complex case, there exists an almost regular Poisson structure $\pi^!$ on $\mathrm{NBl}(M,\mathcal{F}_\pi)$ such that  \begin{enumerate}\item $\mathcal{F}_{\pi^!}=\mathrm{NBl}(\mathcal{F}_\pi)$,\item and that makes $p\colon \left(\mathrm{NBl}(M,\mathcal{F}_\pi), \pi^!\right)\longrightarrow (M, \pi)$ a map of foliated bivector fields. In particular, $p$ is a Poisson map when $\pi$ is a Poisson structure.\end{enumerate}\end{enumerate}\end{proposition}\begin{proof}\begin{enumerate}\item \item Let $\mathbb{K}=\mathbb C$. By compacity of the fibers of $p$ one has $p^*C^{\infty}(M)\equiv C^\infty(\mathrm{NBl}(M, \mathcal{F}))$. This allows to define a Poisson bivector field on $\mathrm{NBl}(M, \mathcal{F})$ in an obvious way  by the formula  \begin{equation}\pi^![p^*f, p^*g]:=p^*\pi[f,g]\; \text{for all}\; f, g\in \mathcal O_M.\end{equation}The bivector field $\pi^!$ clearly satisfies  $(a)$ and $(b)$.\end{enumerate}\end{proof}

\begin{example}[Concentric spheres]\label{ex:spheres}
    Consider the Lie–Poisson structure on $M=~{\mathfrak{so}}(3,\mathbb R)^*$ given on a dual basis $x, y, z$ by

    $$\pi= -z\frac{\partial}{\partial x}\wedge \frac{\partial}{\partial y} + y\frac{\partial}{\partial x}\wedge \frac{\partial}{\partial 
    z} -x\frac{\partial}{\partial y}\wedge \frac{\partial}{\partial z}.$$  The singular foliation $\mathcal{F}_\pi$ is generated by the vector fields 
 $$X:=\pi^\#(dx)=z \frac{\partial}{\partial y}- y\frac{\partial}{\partial z},\;\; Y:=\pi^\#(dy)=z\frac{\partial}{\partial x} - x\frac{\partial}{\partial z},\;\;Z:=\pi^\#(dz)=y \frac{\partial}{\partial x} - x\frac{\partial}{\partial y} .$$
The open dense subset of regular points is $M_{\mathrm{reg}}=\mathbb R^3\setminus \left\{0\right\}$. The kernel of $\pi^{\#}$ at a point $(x,y,z)\in M_{\mathrm{reg}}$ is the line generated by $(x,y,z)$. %A direct computation shows that $$p^{-1}({M}_{\mathrm{reg}})= \left\{ (m,\ell)\in \mathbb R^3\setminus \left\{0\right\}\times \mathbb P^2\,|\, m\in \ell\right\}.$$ 
Hence, the blow-up space $\mathrm{Nbl}_{\mathcal{F}_\pi}(M)=\overline{\left\{ (m,\ell)\in \mathbb R^3\setminus \left\{0\right\}\times \mathbb P^2\,|\, m\in \ell\right\}}$ is the usual blow-up $\mathrm{Bl}_0(\mathbb R^3)$ of $\mathbb R^3$ at $0$. In particular, $ \mathrm{Nbl}_{\mathcal{F}}(M)$ is a smooth manifold.

Let us look at the blow-up space for example in the $x$-chart given by $(x,y,z)\mapsto \left((x, xy, xz),[1: y,z]\right)$ so that $p((x, xy, xz),[1: y,z])=(x, xy, xz)$. The singular foliation $p^!\mathcal{F}_\pi$ is generated by the vector fields  $$\overline{X}=z\frac{\partial}{\partial y}- y\frac{\partial}{\partial z}, \quad\overline{Y}=xz\frac{\partial}{\partial x}-yz\frac{\partial}{\partial y}- (z^2+1)\frac{\partial}{\partial z}$$$\text{and}\;\; \overline{Z}=xy\frac{\partial}{\partial x}-(y^2+1)\frac{\partial}{\partial y}- yz\frac{\partial}{\partial z}$. Since $\overline{X}=y\overline{Y}-z\overline{Z}$, it follows that $p^!\mathcal{F}_\pi$  is a regular foliation on $\mathrm{Nbl}_{\mathcal {F}}(M)$ of rank $2$. We have a short exact sequence of Lie algebroids $$0\rightarrow O(-1)_{\mathrm{Nbl}_\mathcal{F}(M)}\rightarrow p^!\left(\mathfrak {so}(3,\mathbb R)\ltimes \mathbb R^3\right)\rightarrow  p^!\left(\mathfrak {so}(3,\mathbb R)\ltimes \mathbb R^3\right)/O(-1)_{\mathrm{Nbl}_\mathcal{F}(M)} \rightarrow 0$$
    where $O(-1)_{\mathrm{Nbl}_\mathcal{F}(M)}$ denotes the Lie algebra bundle obtained by the restriction of the tautological line bundle $O(-1)\rightarrow \mathbb{P}^{2}\times \mathbb R^3$ to $\mathrm{Nbl}_\mathcal{F}(M)$ and $$p^!\mathcal{F}_\pi\simeq p^!\left(\mathfrak {so}(3,\mathbb R)\ltimes \mathbb R^3\right)/O(-1)_{\mathrm{Nbl}_\mathcal{F}(M)}.$$

The bivector field $\pi$ does not lift to $\mathrm{Nbl}_{\mathcal {F}}(M)$. Recall that the blow-up space of $\mathbb R^3$ at the origin   $\mathrm{Bl}_0(\mathbb R^3)\subset~ \mathbb P^2\times\mathbb R^3$ is covered by three charts given by $x\neq 0$, $y\neq 0$ and $z\neq 0$.  Let us look at the $x$-chart where the projection $p$ becomes $(x,y,z)\mapsto~(x,xy, xz)$. In this chart, $\pi$ pulls back to
   \begin{equation}\label{eq:counter-poisson}
       y\frac{\partial}{\partial z}\wedge \frac{\partial}{\partial x}+z\frac{\partial}{\partial x}\wedge \frac{\partial}{\partial y}+\frac{1}{x}(1+y^2+z^2)\frac{\partial}{\partial y}\wedge \frac{\partial}{\partial x}.
   \end{equation} For $x= 0$, Equation \eqref{eq:counter-poisson} is not defined. In conclusion, the Hamiltonian vector fields of the Poisson structure $(M,\pi)$ lift to $\mathrm{Nbl}_{\mathcal {F}}(M)$ which is a smooth manifold, but the bivector field $\pi$ does not lift to a smooth bivector field on  $\mathrm{Nbl}_{\mathcal {F}}(M)$.

%Notice that the pullback  $\pi^!$ of the bivector field $\pi$  is only defined on the regular part of $\mathrm{Nbl}_\mathcal{F}(M)$ on which $p$ is a diffeomorphism and not on the exceptional divisors $p^{-1}(0)$, for instance in the $x$-chart it is not defined for $x=0$ \cite{Ruben3}. %However, if we replace $\mathbb R$ by $\mathbb C$, $\mathrm{Nbl}(\mathcal{F})$xxxx
\end{example}
In fact, the Poisson structure of  Example \ref{ex:spheres} can be seen as a particular case of a more general construction of Poisson structures in dimension $3$.

\begin{example}\label{ex:non-smooth}
 Let $\varphi\in \mathbb K[x,y,z]$. Consider the Poisson structure on $\mathbb K^3$ given by:

$$\pi=\frac{\partial \varphi}{\partial z}\frac{\partial}{\partial x}\wedge \frac{\partial}{\partial y}-\frac{\partial \varphi}{\partial y}\frac{\partial}{\partial x}\wedge \frac{\partial}{\partial z}+\frac{\partial\varphi}{\partial x}\frac{\partial}{\partial y}\wedge\frac{\partial}{\partial z}$$
This Poisson structure was studied in \cite{PICHEREAU}. Consider $\mathcal{F}_\pi$ its the singular foliation generated by the vector fields $X_{ij}=\frac{\partial \varphi}{\partial x_i}\frac{\partial }{\partial x_j}-\frac{\partial \varphi}{\partial x_j}\frac{\partial }{\partial x_i}$ with $x_i, x_j\in \{x, y, z\}$. When $\frac{\partial\varphi}{\partial x}$, $\frac{\partial \varphi}{\partial y}$ and $\frac{\partial \varphi}{\partial z}$ form a complete intersection, $\mathcal{F}_\pi$ is exactly the space of all $X\in \mathfrak X(\mathbb K^3)$ satisfying $X[\varphi]=0$, see \cite{LLS} Section 3.3. The set of the singular points is the zero locus $W$ of the ideal generated by $\frac{\partial\varphi}{\partial x}$, $\frac{\partial \varphi}{\partial y}$ and $\frac{\partial \varphi}{\partial z}$. For $m=(x,y,z)\in M_{\mathrm{reg}}$, a direct computation shows that $\ker\pi^\#|_{m}=\left[\frac{\partial\varphi}{\partial x}(m):\frac{\partial\varphi}{\partial y}(m):\frac{\partial\varphi}{\partial z}(m)\right]$. Therefore, $\mathrm{Nbl}_{\mathcal{F}_\pi}(M)$ is the natural transformation of $\mathbb K^3$ along $W$ which is not smooth in general. Take, e.g., $\varphi=xy- \frac{z^{n+1}}{n+1}$ for some integer $n\geq 2$: if we denote by $[u:v:w]$ the coordinates of the projective space $\mathbb P^2$,  $\mathrm{Nbl}_{\mathcal{F}_\pi}(M)=\{uy-vx=0, wx-uz^n=0, wy-vz^n=0\}$ is the blow-up of $\mathbb K^3$ along the ideal $(x, y, z^n)$.  In the chart $u=1$, $\mathrm{Nbl}_{\mathcal{F}_\pi}(M)= A^n: wx-z^n=0$ is a Du Val singularity.

By Theorem \ref{thm:main}, even when  $\mathrm{Nbl}_{\mathcal{F}_\pi}(M)$ is singular, it admits a Debord singular foliation $p^!\mathcal{F}_\pi$ whose leaves are manifolds included in $\mathrm{Nbl}_{\mathcal{F}_\pi}(M)$.
\end{example}

We end this section with an example where the blow-up space $\mathrm{Nbl}_\mathcal{F}(M)$ arises as a quotient space. 
\begin{example}[The adjoint action of $\mathfrak {su}(n)$]\label{ex:su(n)}
  Consider the adjoint action of the Lie algebra $\mathfrak {su}(n)$ of traceless anti‑Hermitian $n\times n$ complex matrices. These matrices are diagonalizable in a unitary basis. $M_{\mathrm{reg}}\subset \mathfrak {su}(n)$ is the open dense subset  of matrices whose centralizers have  minimal dimension, i.e., traceless anti‑Hermitian matrices that admit distinct eigenvalues.  %the matrices   in that case the characteristic polynomial and the minimal polynomial coincide.
  By definition $\mathrm{Nbl}_\mathcal F(M)=\overline{\{ C(\mathfrak M), \mathfrak M\in M_{\mathrm{reg}}\subset\mathfrak {su}(n)\}}$, here $C(\mathfrak M)$ is the centralizer of $\mathfrak M\in \mathfrak {su}(n)$.

  Denote by $\mathfrak D\subseteq \mathfrak {su}(n)$ the space of diagonal  anti-Hermitian matrices. We have a well-defined surjective map 
  \begin{equation*}
      \frac{{SU}(n)\times \mathfrak D}{\mathrm{Stab}{(\mathfrak D})}\to \mathfrak {su}(n),\; [(U,D)]\mapsto U^{-1}DU,
  \end{equation*}
  which is proper and invertible on $M_{\mathrm{reg}}$. Notice that for all $\mathfrak M\in M_{\mathrm{reg}}$ such that $\mathfrak M=U^{-1}DU$ for some $(U,D)\in {SU}(n)\times \mathfrak D$, we have $C(U^{-1}DU)=U^{-1}\mathfrak D U$. This implies that there is a well-defined map at the level of the Grassmanian
  \begin{equation}\label{eq:su(n)}
      \frac{{SU}(n)\times\mathfrak D}{\mathrm{Stab}{(\mathfrak D})}\to\mathrm{Nbl}_\mathcal{F}(M),\;[(U,D)]\mapsto (U^{-1} D U,\, U^{-1}\mathfrak D U).
  \end{equation}
Here, $SU(n)$  stands for the Lie group of $n\times n$ unitary matrices with determinant 1, and $\mathrm{Stab}{(\mathfrak D})=\{U\in {SU}(n) \mid U^{-1}\mathfrak D U= \mathfrak D\}$ is the stabilizator of $\mathfrak D$ in $\mathfrak {su}(n)$. The result is similar for $\mathfrak{so}(n)$. Let us show that the map \eqref{eq:su(n)} is a bijection: 

\begin{enumerate}
    \item Surjectivity: let $V\in p^{-1}(\mathfrak M)$ and a convergent sequence of regular points $\mathfrak M_n\underset{n \to +\infty}{\longrightarrow} \mathfrak M$ such that $C(\mathfrak M_n)\underset{n \to +\infty}{\longrightarrow} V$. We have that $\mathfrak M_n=U_n^{-1}D_nU_n\underset{n \to +\infty}{\longrightarrow}~\mathfrak M$ with $(U_n,D_n)\in {SU}(n)\times \mathfrak D$. By compacity of $SU(n)$, there exists a subsequence of  $U_n^{-1}\mathfrak DU_n$ that converges to $V$. This proves surjectivity.

    \item Injectivity: let $U_1, U_2\in SU(n)$ be such that $U_1^{-1}\mathfrak DU_1=U_2^{-1}\mathfrak DU_2$. This implies that $(U_1U_2^{-1})^{-1}\mathfrak D (U_1U_2^{-1})=\mathfrak D$. Thus, $U_1U_2^{-1}\in \mathrm{Stab}{(\mathfrak D})$.
\end{enumerate}
\end{example}

A more general result than Example \ref{ex:su(n)} is the topic of a forth-coming paper.
\section{Conclusion}
We end the paper with a natural question. This paper highlights for a given Lie algebroid $(A, [\cdot,\,\cdot]_A,\rho_A)$ a procedure for separating the space of objects (i.e., the singular foliation part) and the symmetry part (a Lie algebra bundle).  That is, $A$ admits a pullback $p^!A$ along the blow-up space $p\colon \mathrm{Nbl}_\mathcal{F}(M)\rightarrow M$  that comes with a short exact sequence of Lie algebroids \begin{equation}
    \label{eq:extension}0\rightarrow K\rightarrow p^!A\rightarrow \mathcal{D}_\mathcal{F}\rightarrow 0
\end{equation}with $K$ a totally intransitive Lie algebroid (i.e., a  Lie algebra bundle) and $\mathcal{D}_\mathcal{F}$  a Lie algebroid whose anchor map is injective on an open dense subset. % which is essentially the induced singular foliation on the base manifold $\mathrm{Nbl}_\mathcal{F}(M)$.

It is natural to ask whether  $p^!A$ is a semi-direct product $\mathcal{D}_\mathcal{F}\ltimes K$, in other words: 
\begin{question}
    What is the meaning of the curvature $C\colon \wedge^2\mathcal{D}_\mathcal{F}\rightarrow K$,\, $$C(u,v):=s([u,v]_{\mathcal{D}_\mathcal{F}})-[s(u), s(v)]_{p^!A}$$ of a section $s\colon \mathcal{D}_\mathcal{F}\rightarrow p^!A$ of $p^!A\rightarrow \mathcal{D}_\mathcal{F}$ %which is a Lie algebra morphism, i.e., such that $s([u,v]_{\mathcal{D}_\mathcal{F}})=[s(u), s(v)]_{p^!A}$ for all $a,b\in \Gamma(\mathcal{D}_\mathcal{F})$
    ? 
\end{question}

Notice that any section $s\colon \mathcal{D}_\mathcal{F}\rightarrow p^!A$ of \eqref{eq:extension} gives rise to a Lie algebra morphism $$\mathcal{D}_\mathcal{F}\rightarrow \mathrm{Der}(Z(K)),\; \alpha\mapsto\nabla_\alpha\beta:=[s(\alpha), \beta]_{p^!A}$$
i.e., a representation of $\mathcal{D}_\mathcal{F}$ on the
center $Z(K)$ of $K$ \cite{V.Rubtsov}. Moreover, if the isotropy Lie algebras $\{\ker\rho_x\}_{x\in M}$ of  $A$ are Abelian, then  $K$ is Abelian and $(\mathcal{D}_\mathcal{F}, \nabla)$ is a representation of $\mathcal{D}_\mathcal{F}$ on $K$.

\bibliographystyle{alpha}
\bibliography{Nash}
\vfill
\begin{center}
    \textsc{Department of Mathematics, Jilin University, Changchun 130012, Jilin, China, and Mathematics Institute, Georg-August-University of Göttingen, Bunsenstrasse 3-5, Göttingen
37073, Germany}
\end{center}
\end{document}